\documentclass[
reprint,
groupedaddress,
showpacs,preprintnumbers,
amsmath,amssymb,
aps,
pre,
]{revtex4-1}
\usepackage{graphicx}% Include figure files
\usepackage{dcolumn}% Align table columns on decimal point
\usepackage{bm}% bold math
\usepackage{subfigure}
\def \t{\theta}
\def \w{\omega}

\begin{document}
%\preprint{APS/123-QED}
\title{Optimal Design of Minimum-Power Stimuli for Spiking Neurons}% Force line breaks with \\
%\thanks{This work was supported by the NSF under Career Award \#0747877}%
\author{Isuru Dasanayake}
 %\altaffiliation[Also at ]{Physics Department, XYZ University.}%Lines break automatically or can be forced with \\
\email{dasanayakei@seas.wustl.edu}
\author{Jr-Shin Li}%
\email{jsli@seas.wustl.edu}
\affiliation{Department of Electrical and Systems Engineering, Washington University in St.Louis, MO 63130, USA.}
%\collaboration{MUSO Collaboration}%\noaffiliation
%\date{\today}% It is always \today, today,
             %  but any date may be explicitly specified

%%%%%%%%%%%%%%%%%%%%%%%%%%%%%%%%%%%%%%%%%%%%%%%%%%%%%%%%%%%%%%%%%%%%%%%%%%%%%%%%%%%%%%%%%%%%%%%%%%%%%%%%%%%%%%
\begin{abstract}
In this article, we study optimal control problems of spiking neurons whose dynamics are described by a phase model. We design minimum-power current stimuli (controls) that lead to targeted spiking times of neurons, where the cases with unbounded and bounded control amplitude are considered. We show that theoretically the spiking period of a neuron, modeled by phase dynamics, can be arbitrarily altered by a smooth control. However, if the control amplitude is bounded, the range of possible spiking times is constrained and determined by the bound, and feasible spiking times are optimally achieved by piecewise continuous controls. We present analytic expressions of these minimum-power stimuli for spiking neurons and illustrate the optimal solutions with numerical simulations.
\end{abstract}
%%%%%%%%%%%%%%%%%%%%%%%%%%%%%%%%%%%%%%%%%%%%%%%%%%%%%%%%%%%%%%%%%%%%%%%%%%%%%%%%%%%%%%%%%%%%%%%%%%%%%%%%%%%%%%%
% PACS, the Physics and Astronomy
                             % Classification Scheme.
%\keywords{Suggested keywords} maximum principle, spiking neurons, phase models                              %display desired
\maketitle
%\tableofcontents

%%%%%%%%%%%%%%%%%%%%%%%%%%%%%%%%%%%%%%%%%%%%%%%%%%%%%%%%%%%%%%%%%%%%%%%%%%%%%%%%%%%%%%%%%%%%%%%%%%%%%%%%%%%%%%%
\section{Introduction}
\label{sec:intro}
Control of neurons and hence the nervous system  by external current stimuli (controls) has received increased scientific attention in recent years for its wide range of applications from deep brain stimulation to oscillatory neurocomputers \cite{ uhlhaas06,osipov07,Izhikevich99}. Conventionally, neuron oscillators are represented by phase-reduced models, which form a standard nonlinear system \cite{Brown04,Winfree01}. Intensive studies using phase models have been carried out, for example, on the investigation of the patterns of synchrony that result from the type and architecture of coupling \cite{Ashwin92, Taylor98} and on the response of large groups of oscillators to external stimuli \cite{Moehlis06, Tass89}, where the inputs to the neuron systems were initially defined and the dynamics of neural populations were analyzed in detail.

Recently, control theoretic approaches have been employed to design external stimuli that drive neurons to behave in a desired way. For example, a multilinear feedback control technique has been used to control the individual phase relation between coupled oscillators \cite{kano10}; a nonlinear feedback approach has been employed to engineer complex dynamic structures and synthesize delicate synchronization features of nonlinear systems \cite{Kiss07}; and our recent work has illustrated controllability of a network of neurons with different natural oscillation frequencies adopting tools from geometric control theory \cite{Li_NOLCOS10}.

There has been an increase in the demand for controlling not only the collective behavior of a network of oscillators but also the behavior of each individual oscillator. It is feasible to change the spiking periods of oscillators or tune the individual phase relationship between coupled oscillators by the use of electric stimuli \cite{Schiff94, kano10}. Minimum-power stimuli that elicit spikes of a neuron at specified times close to the natural spiking time were analyzed \cite{Moehlis06}. Optimal waveforms for the entrainment of weakly forced oscillators that maximize the locking range have been calculated, where first and second harmonics were used to approximate the phase response curve \cite{kiss10}. These optimal controls were found mainly based on the calculus of variations, which restricts the optimal solutions to the class of smooth controls and the bound of the control amplitude was not taken into account.

In this paper, we apply the Pontryagin's maximum principle \cite{Pontryagin62, Stefanatos10} to derive minimum-power controls that spike a neuron at desired time instants. We consider both cases when the available control amplitude is unbounded and bounded. The latter is of practical importance due to physical limitations of experimental equipment and the safety margin for neurons, e.g., the requirement of a mild brain stimulations in neurological treatments for Parkinson's disease and epilepsy.

This paper is organized as follows. In Section \ref{sec:phase_model}, we introduce the phase model for spiking neurons and formulate the related optimal control problem. In Section \ref{sec:minpower_control}, we derive minimum-power controls associated with specified spiking times in the absence and presence of control amplitude constraints, in which various phase models including sinusoidal PRC, SNIPER PRC, and theta neuron models are considered. In addition, we present examples and simulations to demonstrate the resulting optimal control strategies.

%%%%%%%%%%%%%%%%%%%%%%%%%%%%%%%%%%%%%%%%%%%%%%%%%%%%%%%
\section {Optimal Control of Spiking Neurons}
\label{sec:phase_model}
A periodically spiking or firing neuron can be considered as a periodic oscillator governed by the nonlinear dynamical equation of the form
\begin{equation}
    \label{eq:phasemodel}
    \frac{d\theta}{dt}=f(\theta)+Z(\theta)I(t),
\end{equation}
where $\theta$ is the phase of the oscillation, $f(\theta)$ and $Z(\theta)$ are real-valued functions giving the neuron's baseline dynamics and its phase response, respectively, and $I(t)$ is an external current stimulus \cite{Brown04}. The nonlinear dynamical system described in \eqref{eq:phasemodel} is referred to as the phase model for the neuron. The assumption that $Z(\theta)$ vanishes only on isolated points and that $f\left(\theta\right)>0$ are made so that a full revolution of the phase is possible. By convention, neuron spikes occur when $\theta=2n\pi$, where $n\in\mathbb{N}$. In the absence of any input $I(t)$, the neuron spikes periodically at its natural frequency, while the periodicity can be altered in a desired manner by an appropriate choice of $I(t)$.

In this article, we study optimal design of neural inputs that lead to the spiking of neurons at a specified time $T$ after spiking at time $t=0$. In particular, we find the stimulus that fires a neuron with minimum power, which is formulated as the following optimal control problem,
\begin{align}
    \label{eq:opt_con_pro}
    \min_{I(t)} \quad & \int_0^T I(t)^2\,dt\\
    {\rm s.t.} \quad & \dot{\theta}=f(\theta)+Z(\theta)I(t), \nonumber\\
    &\theta(0)=0, \quad \theta(T)=2\pi \nonumber\\
    &|I(t)|\leq M, \ \ \forall\ t, \nonumber
\end{align}
where $M>0$ is the amplitude bound of the current stimulus $I(t)$. Note that instantaneous or arbitrarily delayed spiking of a neuron is possible if $I(t)$ is unbounded, i.e., $M=\infty$; however, the range of feasible spiking periods of a neuron described as in \eqref{eq:phasemodel} is restricted with a finite $M$. We consider both unbounded and bounded cases.

%%%%%%%%%%%%%%%%%%%%%%%%%%%%%%%%%%%%%%%%%%%%%%%%%%%%%%%%%%%%%%%%%%%%%%%%%%%%%%%%%%%%%%%%%%%%%%%%%%%%%%%%%%%%%%%%%%%
\section{Minimum-Power Stimulus for Specified Firing Time}% with Unbounded Amplitude}
\label{sec:minpower_control}
We consider the minimum-power optimal control problem of spiking neurons as formulated in \eqref{eq:opt_con_pro} for various phase models including sinusoidal PRC, SNIPER PRC, and theta neuron.

\subsection{Sinusoidal PRC Phase Model}
\label{sec:sine_prc}
Consider the sinusoidal PRC model,
\begin{equation}
    \label{eq:sin_model}
    \dot{\theta}=\omega+z_d\sin\theta\cdot I(t),
\end{equation}
where $\omega$ is the natural oscillation frequency of the neuron and $z_d$ is a model-dependent constant. The neuron described by this phase model spikes periodically with the period $T=2\pi/\omega$ in the absence of any external input, i.e., $I(t)=0$.

%%%%%%%%%%%%%%%%%%%%%%%%%%%%%%%%%%%%%%%%%%%%%%%%%%%%%%%%%%%%%%%%%%
\subsubsection{Spiking Neurons with Unbounded Control}
\label{sec:unbounded_control_sine}
The optimal current profile can be derived by Pontryagin's Maximum Principle \cite{Pontryagin62}. Given the optimal control problem as in \eqref{eq:opt_con_pro}, we form the control Hamiltonian
\begin{equation}
    \label{eq:hamiltonian}
	H=I^2+\lambda(\omega+z_d\sin\theta\cdot I),
\end{equation}
where $\lambda$ is the Lagrange multiplier. The necessary optimality conditions according to the Maximum Principle give
\begin{align}
	\label{eq:lambda_dot}
	\dot{\lambda}=-\frac{\partial H}{\partial \theta}=-\lambda z_d I\cos\theta,
\end{align}
and $\frac{\partial H}{\partial I}=2I+\lambda z_d\sin\theta=0$. Hence, the optimal current $I$ satisfies
\begin{eqnarray}
	\label{eq:I}
	I=-\frac{1}{2}\lambda z_d\sin\theta.
\end{eqnarray}

Substituting \eqref{eq:I} into \eqref{eq:sin_model} and \eqref{eq:lambda_dot}, the optimal control problem is then transformed to a boundary value problem, which characterizes the optimal trajectories of $\theta(t)$ and $\lambda(t)$,
\begin{align}
	\label{eq:theta_p2}
    \dot{\theta} &= \omega-\frac{z_d^2\lambda}{2} \sin^2 \theta,\\
	\label{eq:lambda_p2}
    \dot{\lambda} &= \frac{z_d^2\lambda^2}{2} \sin\theta\cos\theta,
\end{align}
with boundary conditions $\theta(0)=0$ and $\theta(T)=2\pi$ while $\lambda(0)$ and $\lambda(T)$ are unspecified.

Additionally, since the Hamiltonian is not explicitly dependent on time, the optimal triple $(\lambda,\theta, I)$ satisfies $H(\lambda,\theta,I)=c$, $\forall\, 0\leq t\leq T$,
where $c$ is a constant. Together with \eqref{eq:I}, this yields
\begin{equation}
	\label{eq:quadratic}
	-\frac{z_d^2}{4}\sin^2\theta\lambda^2+\omega\lambda=c,
\end{equation}
and the constant $c=\omega\lambda_0$ is obtained by the initial conditions $\theta(0)=0$ and $\lambda(0)=\lambda_0$, which is undetermined. Then, the optimal multiplier can be found by solving the above quadratic equation \eqref{eq:quadratic}, which gives
\begin{equation}
	\label{eq:lambda}
\lambda=\frac{2\omega\pm2\sqrt{\omega^2-\omega\lambda_0z_d^2\sin^2\theta}}{z_d^2\sin^2\theta},
\end{equation}
and the optimal trajectory of $\theta$ follows
\begin{equation}
	\label{eq:theta}
	\dot{\theta}=\mp\sqrt{\omega^2-\omega\lambda_0 z_d^2\sin^2 \theta},
\end{equation}
by plugging $\lambda$ in \eqref{eq:lambda} into \eqref{eq:theta_p2}. Integrating \eqref{eq:theta} by separation of variables, we find the spiking time $T$ with respect to the initial condition $\lambda_0$,
\begin{equation}
	\label{eq:T}  T=\frac{1}{\omega}F\left(\scriptstyle{2\pi,\sqrt{\frac{\lambda_0}{\omega}}z_d}\right)=\int_0^{2\pi}{\frac{1}{\scriptstyle{\sqrt{\omega^2-\omega\lambda_0 z_d^2\sin^2\theta}}}}d\theta,
\end{equation}
where $F$ denotes the elliptic integral. Note that we choose the positive sign in \eqref{eq:theta} since the negative velocity indicates the backward phase evolution. Therefore, given a desired spiking time $T$ of the neuron, the initial value, $\lambda_0$, corresponding to the optimal trajectory of the multiplier can be found via the one-to-one relation in \eqref{eq:T}. Consequently, the optimal trajectories of $\theta$ and $\lambda$ can be easily computed by evolving \eqref{eq:theta_p2} and \eqref{eq:lambda_p2} forward in time. Plugging \eqref{eq:lambda} into \eqref{eq:I}, we obtain the optimal feedback law
\begin{equation}
    \label{eq:I*}
    I^*= \frac{-\w+\sqrt{\w^2-\w\lambda_0 z_d^2\sin^2\t}}{z_d\sin\t},
\end{equation}
which drives the neuron from $\theta(0)=0$ to $\theta(T)=2\pi$ with minimum power.

\begin{figure}[t]
\centering
\includegraphics[scale=0.2]{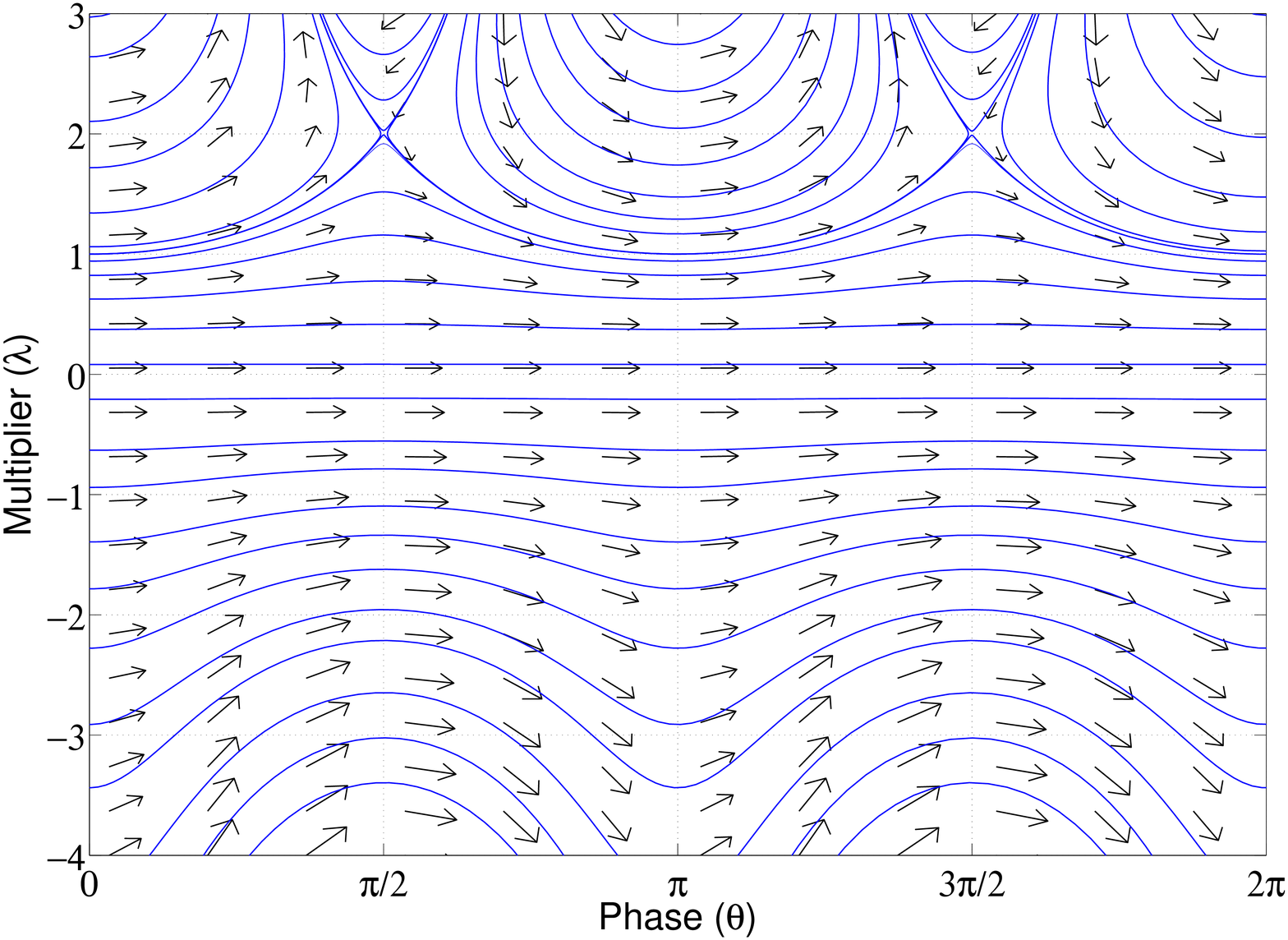}
\caption{Extremals of sinusoidal PRC system with $z_d=1$ and $\w=1$}
\label{fig:vectorfield}
\end{figure}

\begin{figure}[b]
\centering
\includegraphics[scale=0.2]{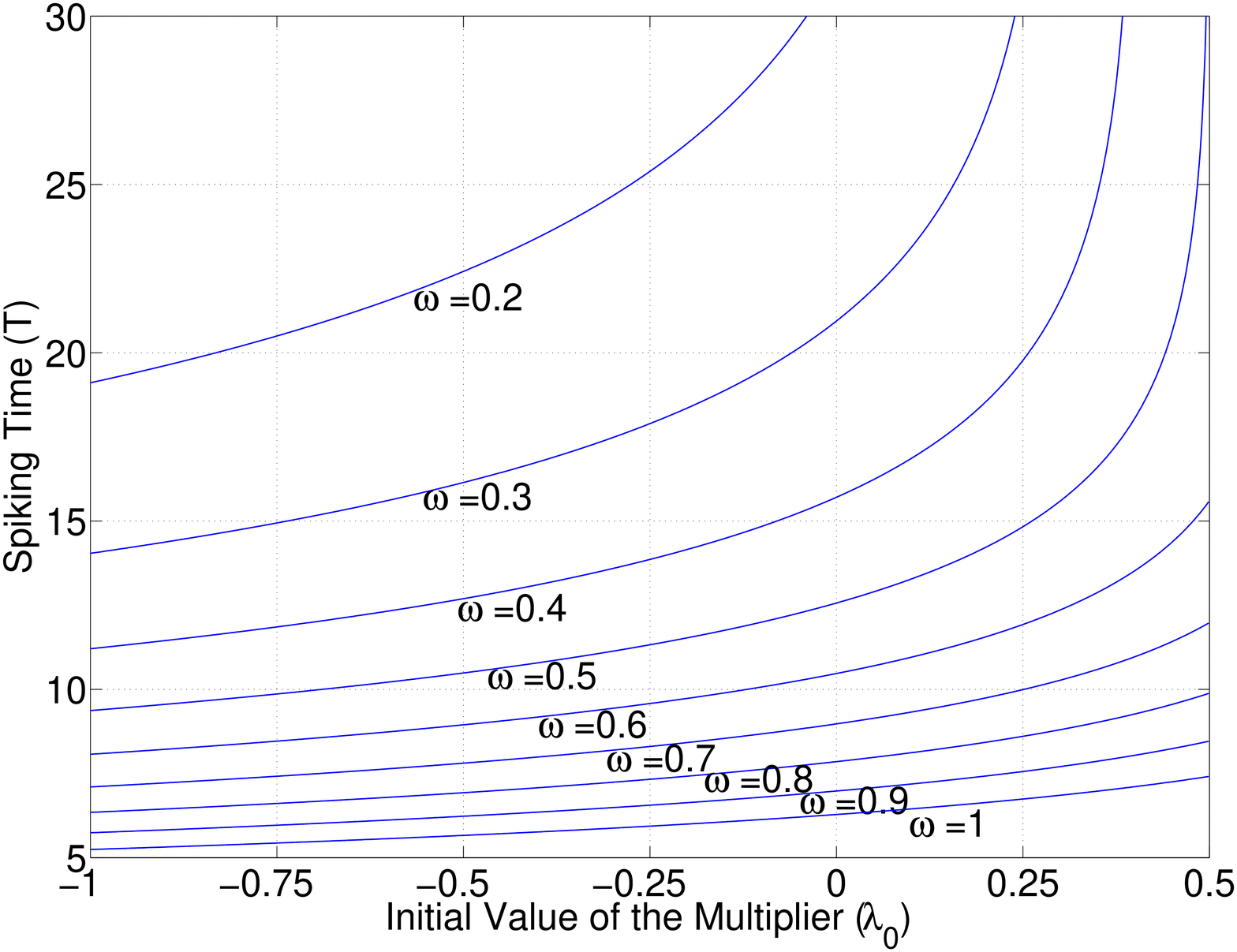}
\caption{Variation of the spiking time, $T$, with respect to the initial multiplier value, $\lambda_0$, leading to optimal trajectories, with different values of $\omega$ and $z_d=1$ for sinusoidal PRC model.}
\label{fig:T_vs_lambda0_sine}
\end{figure}

The feasibility of spiking the neuron at a desired time $T$ largely depends on the initial value of the multiplier, $\lambda_0$. It is clear from \eqref{eq:T} and \eqref{eq:theta} that a complete $2\pi$ revolution is impossible when $\lambda_0>\omega/z_d^2$. This fact also can be seen from FIG. \ref{fig:vectorfield}, where the system evolution defined by \eqref{eq:theta_p2} and \eqref{eq:lambda_p2} for $z_d=1$ and $\omega=1$ with respect to different $\lambda_0$ values ($\t=0$ axis) is illustrated. When $\lambda_0=0$, according to equation \eqref{eq:T}, the spiking period is equal to the natural spiking period, $2\pi/\omega$, and no external stimulus needs to be applied, i.e., $I^*(t)=0$, $\forall t\in[0,2\pi/\w]$. Since, from \eqref{eq:T}, $T$ is a monotonically increasing function of $\lambda_0$ for fixed $\w$ and $z_d$ and, from \eqref{eq:theta}, the average phase velocity decreases when $\lambda_0$ increases, the spiking time $T>2\pi/\omega$ for $\lambda_0>0$ and $T<2\pi/\omega$ for $\lambda_0<0$. FIG. \ref{fig:T_vs_lambda0_sine} shows variation of the spiking time $T$ with the $\lambda_0$ corresponding to the optimal trajectories for different $\omega$ values with $z_d=1$.

% ================== Figure 1 ================

% ================== Figure 2 ================

% ================== Figure 3 ================
\begin{figure}[t]
\subfigure[]{\includegraphics[scale=0.105]{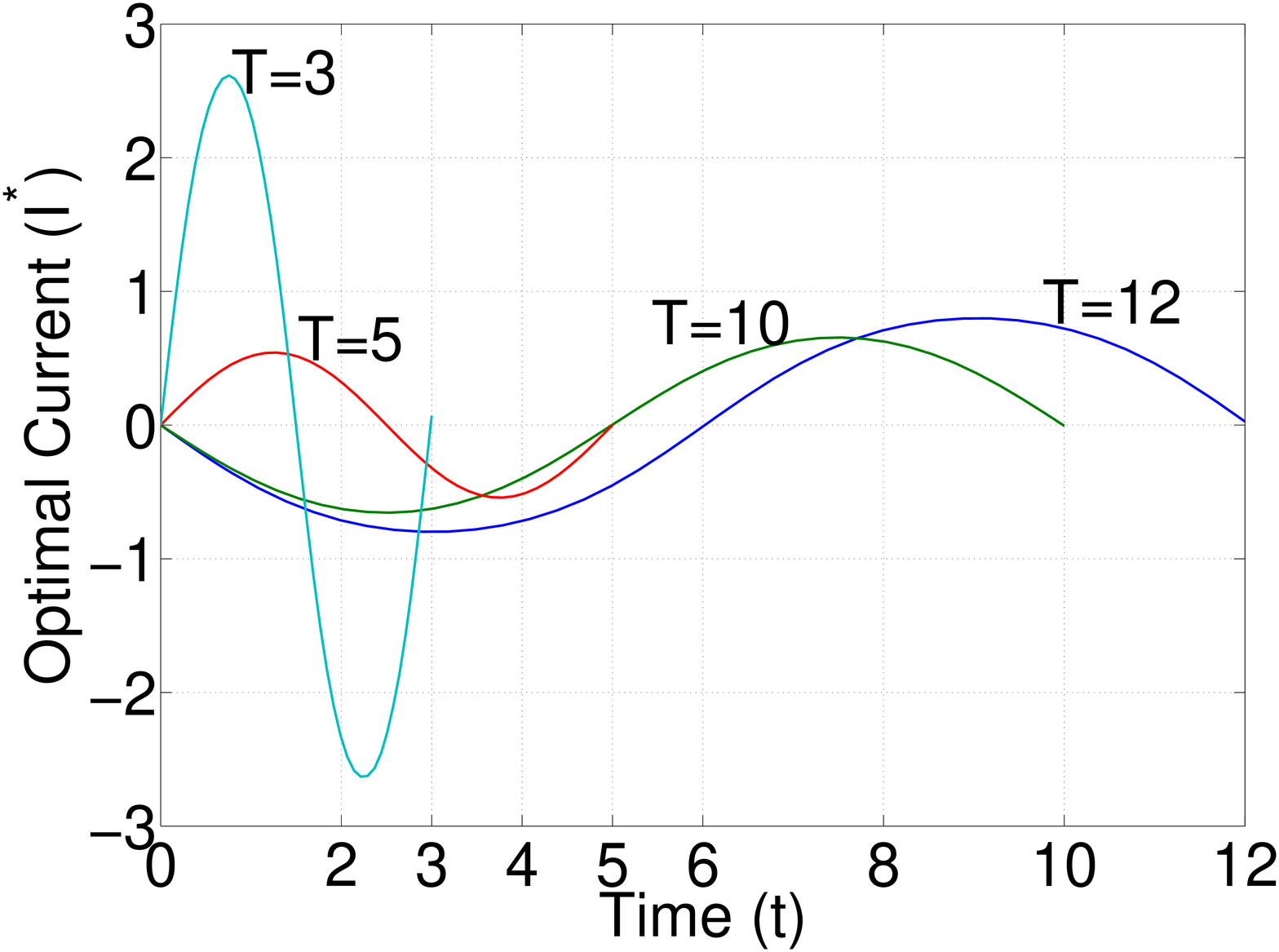}\label{fig:control_vs_time_sine}}
\subfigure[]{\includegraphics[scale=0.105]{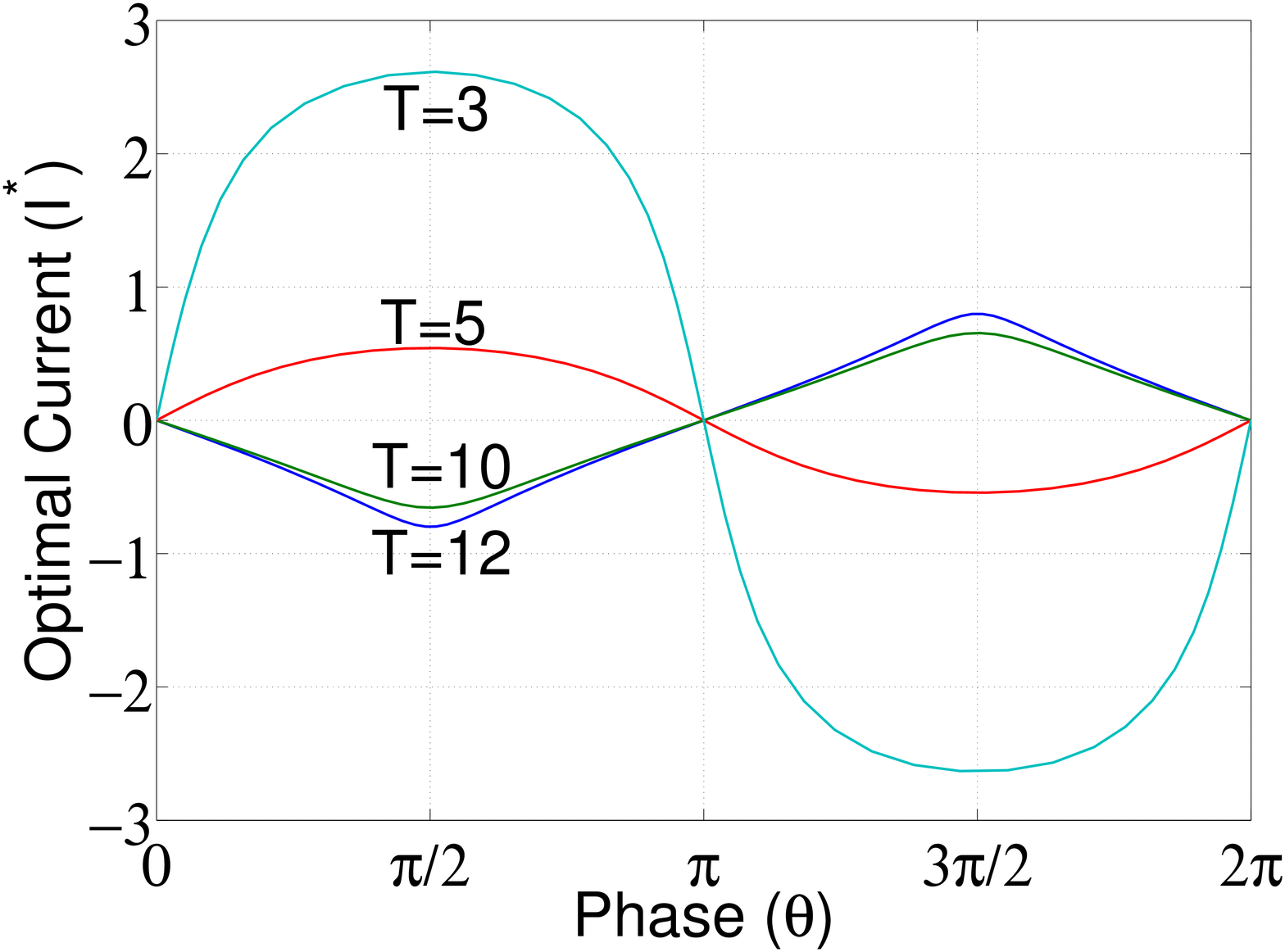}\label{fig:control_vs_theta_sine}}\\ \subfigure[]{\includegraphics[scale=0.105]{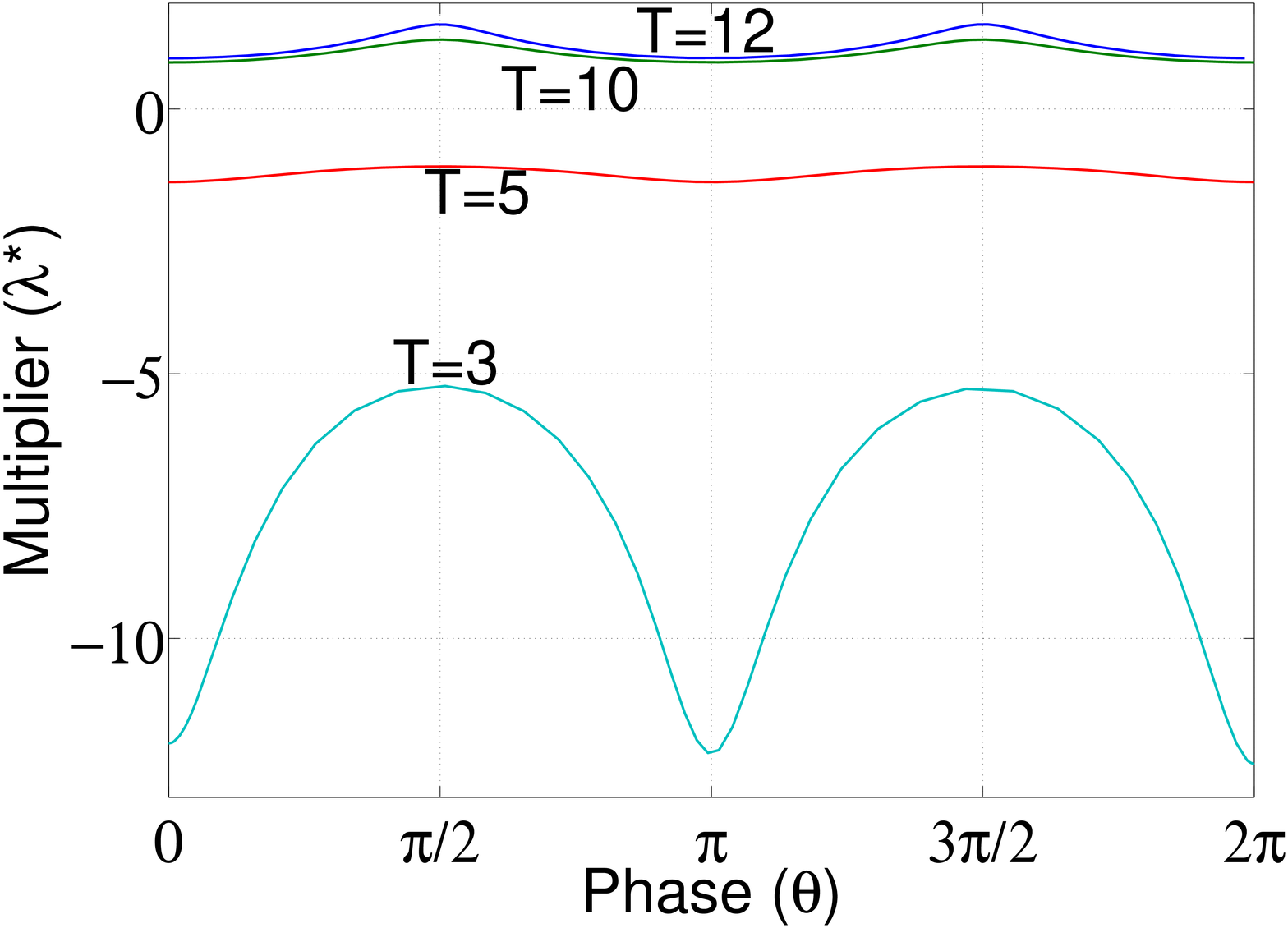}\label{fig:lambda_vs_theta_sine}}
\subfigure[]{\includegraphics[scale=0.105]{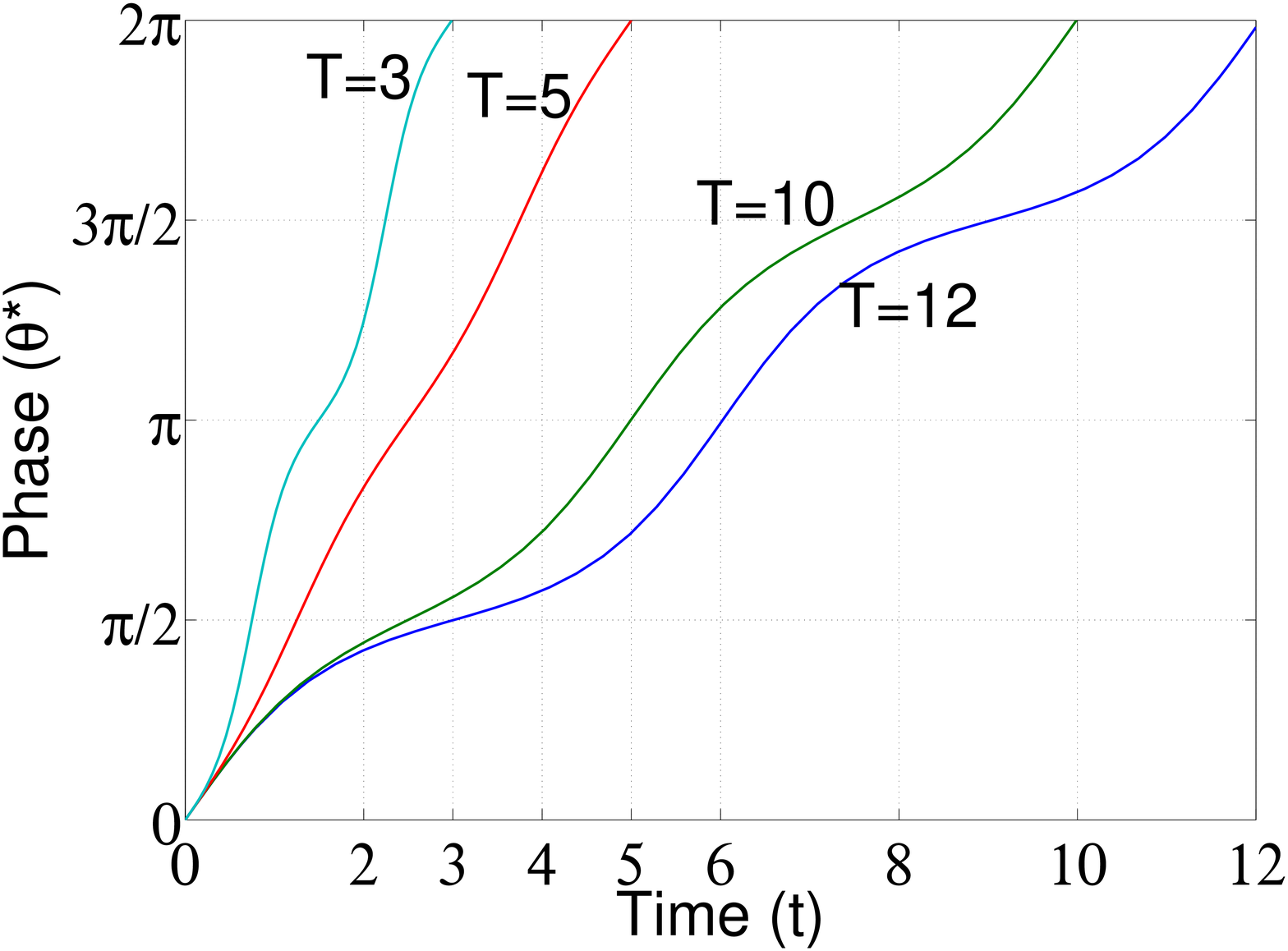}\label{fig:theta_vs_time_sine}}
	  \caption{Optimal solutions for various spiking times $T=3,5,10,12$ for sinusoidal PRC model with $z_d=1$ and $\w=1$. \subref{fig:control_vs_time_sine} The minimum-power control $I^*$. \subref{fig:control_vs_theta_sine} Variation of $I^*$ with phase $\theta$. \subref{fig:lambda_vs_theta_sine} Variation of the optimal multiplier, $\lambda^*$, with $\theta$. \subref{fig:theta_vs_time_sine} Optimal phase trajectories following $I^*$.}
	\label{fig:trajectories_sine}
\end{figure}

%==============================================================================================================================================
The relation between the spiking time $T$ and required minimum power $E=\min\int_0^{T}{{I}^2(t)}dt$ is evident via a simple sensitivity analysis \cite{Bryson75}. Since a small change in the initial condition, $d\theta$, and a small change in the initial time, $dt$, result in a small change in power according to $dE=\lambda(t) d\theta-H(t)d\theta$, it follows that \cite{Bryson75}
\begin{equation}
    \label{eq:energy}
    -\frac{\partial E}{\partial t}=H=c=\w\lambda_0.
\end{equation}
This implies that $E$ increases with initial time $t$ for $\lambda_0<0$ and decreases for $\lambda_0>0$. Since the increment of the initial time is equivalent to the decrement of spiking time $T$, $\partial E/\partial T=\omega\lambda_0$. Since $\lambda_0<0$ ($\lambda_0>0$) corresponds to $T<2\pi/\omega$ ($T>2\pi/\omega$), we see that the required minimum power increases if we move away from the natural spiking time.

The minimum-power stimulus $I^*$ as in \eqref{eq:I*} plotted with respect to time and phase for various spiking times $T=3,5,10,12$ with $\w=1$ and $z_d=1$ are shown in FIG. \ref{fig:control_vs_time_sine} and \ref{fig:control_vs_theta_sine}, respectively. The respective optimal trajectories of $\lambda(\t)$ and $\t(t)$ for these spiking times are presented in FIG. \ref{fig:lambda_vs_theta_sine} and \ref{fig:theta_vs_time_sine}.

%%%%%%%%%%%%%%%%%%%%%%%%%%%%%%%%%%%%%%%%%%%%%%%%%%%%%%%%%%%%%%%%%%%%%%%%%%%%%%%%%%%%%%%%%%%%%%%%%%%%%%%%%%%%
\subsubsection{Spiking Neurons with Bounded Control}
\label{sec:bounded_control_sine}
In practice, the amplitude of stimuli in physical systems are limited, so we consider spiking the sinusoidal neuron with bounded control amplitude, namely, in the optimal control problem \eqref{eq:opt_con_pro}, $|I(t)|\leq M<\infty$ for all $t\in[0,T]$, where $T$ is the desired spiking period. In this case, there exists a range of feasible spiking periods depending on the value of $M$, in contrast to the previous case where any desired spiking time is feasible. We first observe that given this bound $M$, the minimum time it takes to spike a neuron can be achieved by choosing the control that keeps the phase velocity $\dot{\t}$ maximum over $t\in[0,T]$. Such a time-optimal control, for $z_d>0$, can be characterized by a switching, i.e.,
\begin{equation}
    \label{sin_Tmin_control}
    I^*_{Tmin}= \left\{\begin{array}{c}  M \quad \mathrm{for} \quad 0\leq \theta <\pi\\ -M \quad \mathrm{for} \quad \pi \leq \theta <2\pi  \end{array} \right..
\end{equation}
Consequently, the spiking time with $I^*_{Tmin}$ can be computed using \eqref{eq:sin_model} and \eqref{sin_Tmin_control}, which yields
\begin{equation}
    T^{M}_{min}= \frac{2\pi-4\tan^{-1}\left\{z_dM/\sqrt{-z_d^2M^2+\omega^2}\right\}}{\sqrt{-z_d^2M^2+\omega^2}}. \label{eq:sin_min_T}
\end{equation}
It follows that $I^*$, derived in \eqref{eq:I*}, is the minimum-power stimulus that spikes the neuron at a desired spiking time $T$ if $|I^*|\leq M$ for all $t\in [0,T]$. However, there exists a shortest possible spiking time by $I^*$ given the bound $M$. Simple first and second order optimality conditions applied to \eqref{eq:I*} find that the maximum value of $I^*$ occurs at $\theta=\pi/2$ for $\lambda_0<0$ and at $\theta=3\pi/2$ for $\lambda_0>0$. Therefore, the $\lambda_0$ for the shortest spiking time with control $I^*$ satisfying $|I^*(t)|\leq M$ can be calculated by substituting $I^*=M$ and $\theta=\pi/2$ to the equation \eqref{eq:I*}, and then from \eqref{eq:T} we obtain this shortest spiking period
\begin{equation}
    \label{eq:min_time_UBC}
    T^{I^*}_{min}=\int_0^{2\pi}\frac{1}{\sqrt{\omega^2+z_dM(z_dM+2\omega)\sin^2(\theta)}}.
\end{equation}
Note that $T^{M}_{min}<T^{I^*}_{min}$. According to \eqref{eq:sin_model} when $M\geq\omega/z_d$, arbitrarily large spiking times can be achieved by making $\dot{\theta}$ arbitrary close to zero. Therefore we consider two cases for $M\geq\omega/z_d$ and $M<\omega/z_d$.

% =================== Figure 4 ====================
\begin{figure}[t]
	\includegraphics[scale=0.2]{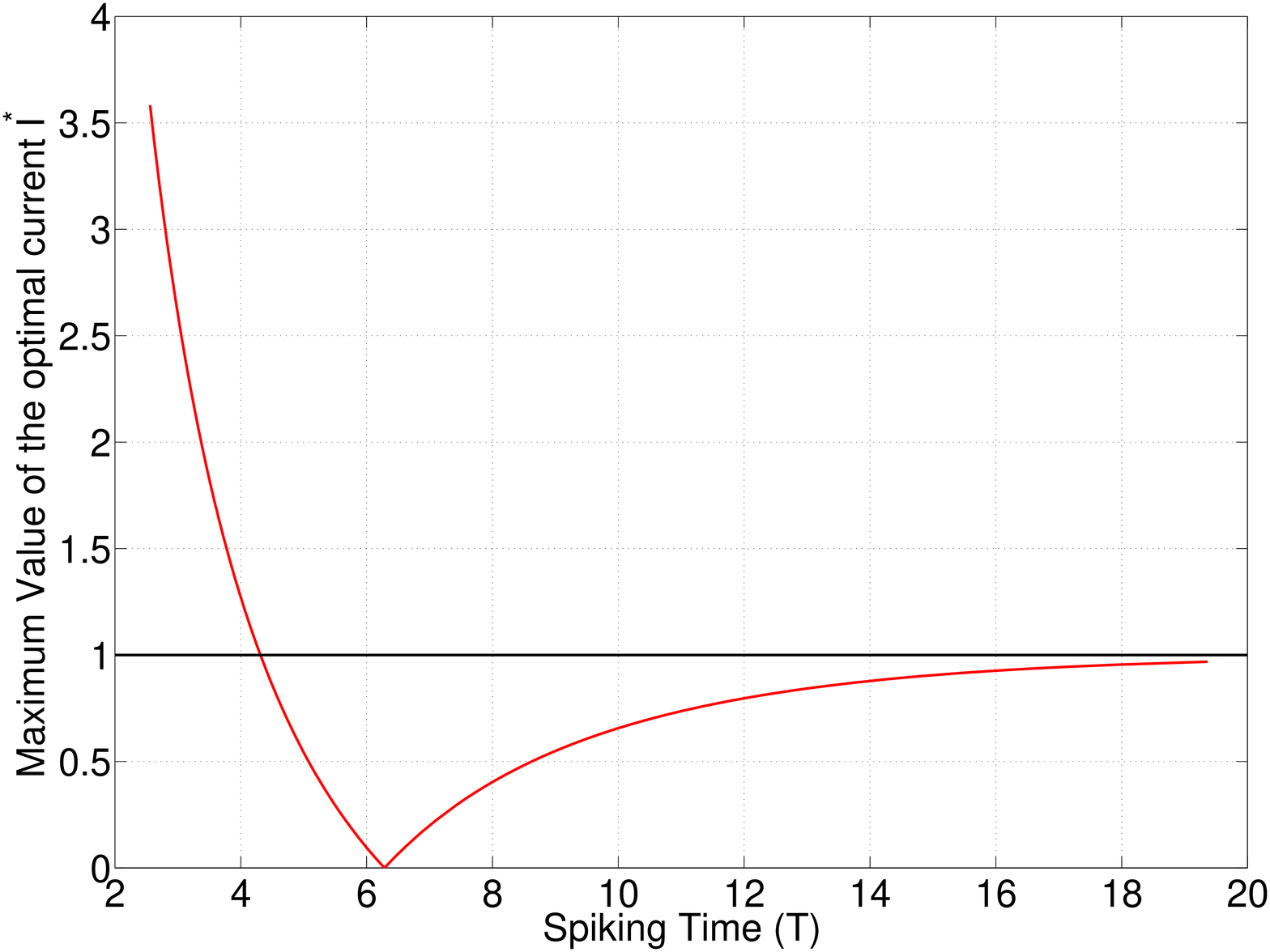}
	\caption{Variation of the maximum value of $I^*$ with spiking time $T$ for sinusoidal PRC model with $\w=1$ and $z_d=1$.}
	\label{fig:max_I_vs_T}
\end{figure}

%----------------------------------------------
\emph{Case I: $M\geq\omega/z_d$}. Since $I^*$ takes the maximum value at $\t=3\pi/2$ for $\lambda_0>0$, we have $|I^*|\leq(\omega-\sqrt{\omega^2-\omega\lambda_0z_d^2})/z_d$, which leads to $|I^*|<\w/z_d\leq M$ for $\lambda_0>0$. This implies that $I^*$ is the minimum-power control for any desired spiking time $T>2\pi/\w$ when $M\geq\w/z_d$, and hence for any spiking time $T\geq T^{I^*}_{min}$. Variation of the maximum value of the control $I^*$ with spiking time $T$ for $\omega=1$ and $z_d=1$ is depicted in FIG. \ref{fig:max_I_vs_T}. Shorter spiking times $T\in[T^{M}_{min}, T^{I^*}_{min})$ are feasible but, due to the bound $M$, can not be achieved by $I^*$ since it requires a control with amplitude greater than $M$ for some $t\in[0,T]$. However, these spiking times can be optimally achieved by applying controls switching between $I^*$ and $I^*_{Tmin}$.

Let the desired spiking time $T\in[T^{M}_{min}, T^{I^*}_{min})$. Then, there exist two angles $\theta_1=\sin^{-1}[-2M\omega/(z_dM^2+z_d\omega \lambda_0)]$ and $\theta_2=\pi-\theta_1$ where $I^*$ meets the bound $M$. When $\theta\in(\theta_1,\theta_2)$, $I^*>M$ and we take $I(\theta)=M$ for $\theta\in[\theta_1,\theta_2]$. The Hamiltonian of the system when $\theta\in[\theta_1,\theta_2]$ is then, from \eqref{eq:hamiltonian},
$H=M^2+\lambda(\omega+z_d\sin\theta\,M)$.
If the triple $(\lambda,\t,M)$ is optimal, then $H$ is a constant, which gives $\lambda=(H-M^2)/(\omega+z_dM\sin\theta)$. This multiplier satisfies the adjoint equation \eqref{eq:lambda_dot}, and therefore $I(\theta)=M$ is optimal for $\theta\in[\theta_1, \theta_2]$. Similarly, by symmetry, $I^*<-M$ when $\t\in[\t_3,\t_4]$, where $\theta_3=\pi+\theta_1$ and $\theta_4=2\pi-\theta_1$, if the desired spiking time $T\in[T^{M}_{min}, T^{I^*}_{min})$. It can be easily shown by the same fashion that $I(\t)=-M$ is optimal in the interval $\t\in[\t_3,\t_4]$.

Therefore, the minimum-power optimal control that spikes the neuron at $T\in[T^{M}_{min}, T^{I^*}_{min})$ can be characterized by four switchings between $I^*$ and $M$, i.e.,
\begin{align}
    \label{eq:I1*}
    I^*_1=\left\{\begin{array}{ll} I^* & \ 0\leq \theta < \theta_1 \\ M & \ \theta_1\leq\theta\leq\theta_2 \\ I^* & \ \theta_2 < \theta < \theta_3 \\ -M & \ \theta_3\leq\theta\leq\theta_4\\ I^* & \ \theta_4 < \theta \leq 2\pi.
    \end{array}\right.
\end{align}
The initial value of the multiplier, $\lambda_0$, resulting in the optimal trajectory, can then be found according to the desired spiking time $T\in[T^{M}_{min}, T^{I^*}_{min})$ through the relation
$$T=\int_0^{\theta_1}{\frac{4}{\sqrt{\scriptstyle{\omega^2-\omega\lambda_0z_d^2\sin^2\theta}}}}d\theta +\int_{\theta_1}^{\frac{\pi}{2}}{\frac{4}{\omega+z_dM\sin\left(\theta\right)}}d\theta.$$

% ====================== Figure 5 ============================================================================
\begin{figure}[t]
\subfigure[]{\includegraphics[scale=0.104]{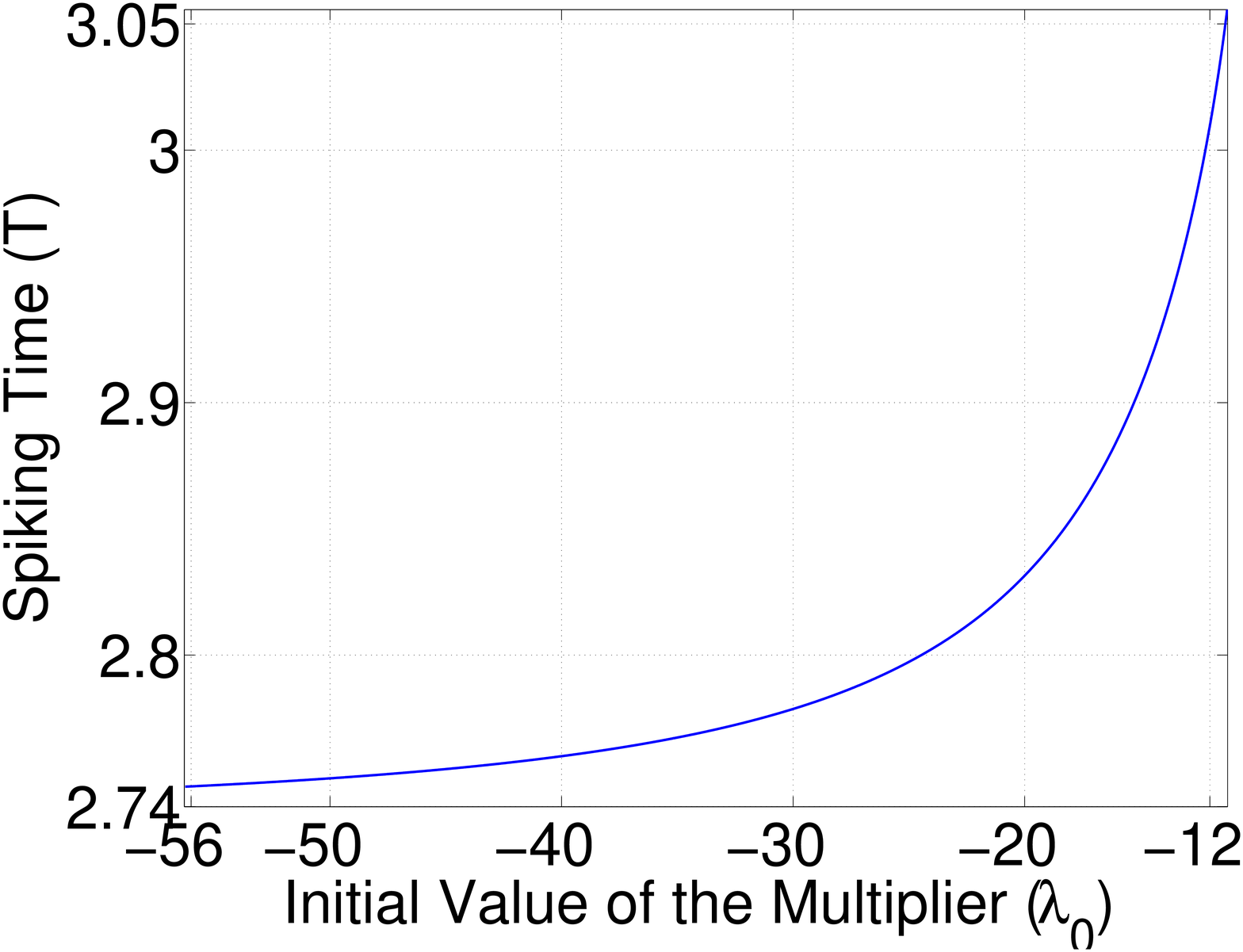}\label{fig:T_vs_lambda0_sine_for_I1}}
\subfigure[]{\includegraphics[scale=0.104]{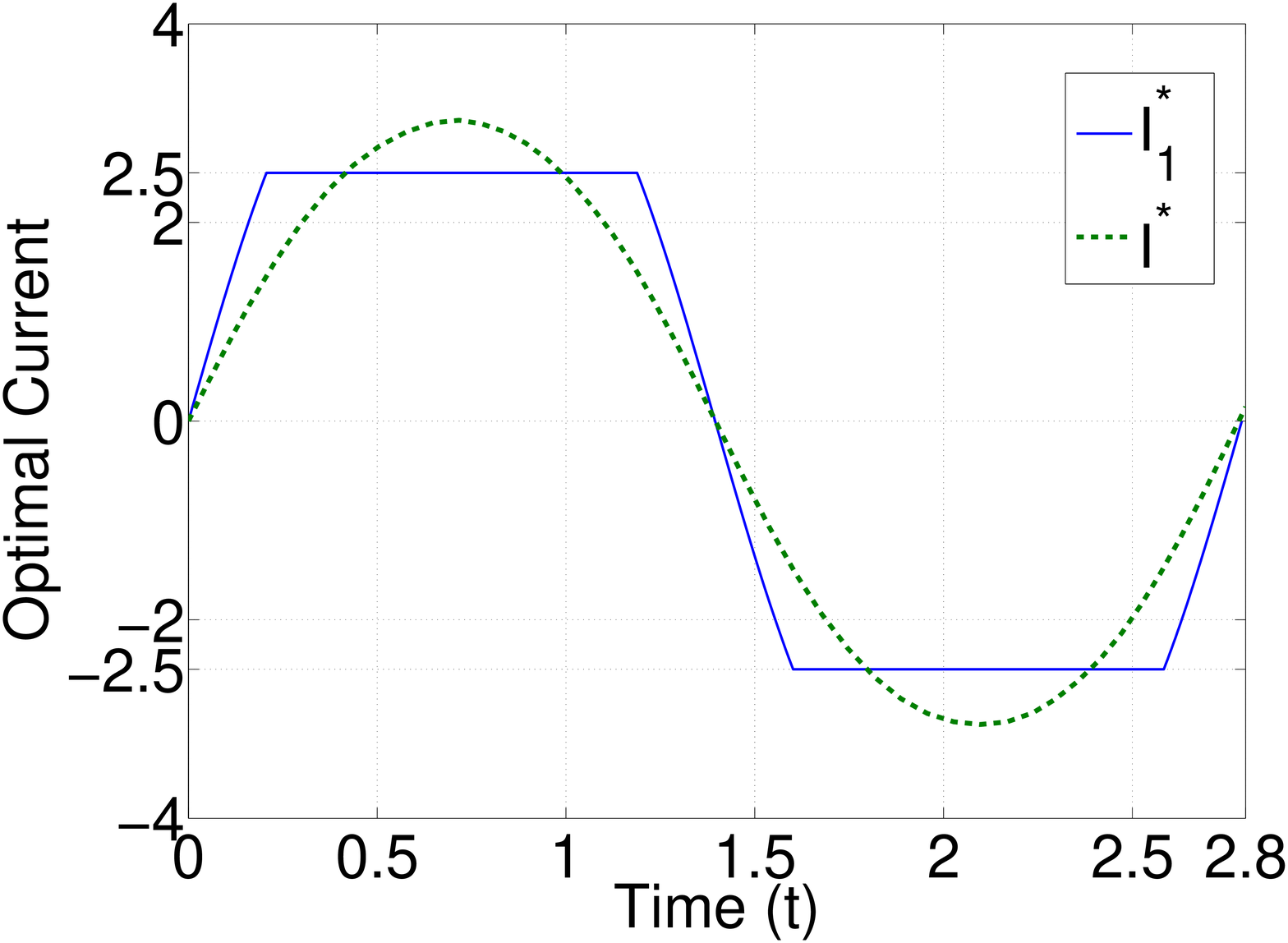}\label{fig:I1_and_I_sine}}
\caption{\subref{fig:T_vs_lambda0_sine_for_I1} Variation of the spiking time $T\in[T^M_{min},T^{I^*}_{min})$ for sinusoidal PRC model with initial multiplier value, $\lambda_0$, for the bound of control amplitude $M=2.5$. \subref{fig:I1_and_I_sine} Minimum-power controls with ($M=2.5$) and without a constraint on the control amplitude for sinusoidal PRC model with $T=2.8$, $z_d=1$, and $\w=1$.}
	\label{fig:bounded_case1_sine}
\end{figure}

FIG. \ref{fig:T_vs_lambda0_sine_for_I1} shows the relation between $\lambda_0$ and $T$ by $I^*_1$ for $M=2.5,\ z_d=1,$ and $ \omega =1$. From \eqref{eq:sin_min_T} the minimum possible spiking time with this control bound $M=2.5$ is $T^{M}_{min}=2.735$ and from \eqref{eq:min_time_UBC} the minimum spiking time by $I^*$ is $T^{I^*}_{min}=3.056$. Thus, in this example, any desired spiking time $T>3.056$ can be optimally achieved by $I^*$ whereas any $T\in[2.735,3.056)$ can be optimally obtained by $I^*_1$ as in \eqref{eq:I1*}. FIG. \ref{fig:I1_and_I_sine} illustrates the bounded and unbounded optimal controls that fire the neuron at $T=2.8$, where $I^*$ is the minimum-power stimulus when the control amplitude is not limited and $I^*_1$ is the minimum-power stimulus when the bound $M=2.5$. $I^*$ drives the neuron from $\t(0)=0$ to $\t(2.8)=2\pi$ with 13.54 units of power whereas $I^*_1$ requires 14.13 units.

%----------------------------------------------
\emph{Case II : $M<\omega/z_d$}. In contrast with Case I in the previous section, achieving arbitrarily large spiking times is not feasible with a bound $M<\omega/z_d$. In this case, the longest possible spiking time is achieved by
$$I^*_{Tmax}= \left\{\begin{array}{ll}  -M & \quad \text{for} \ \ 0\leq \theta <\pi, \\ M & \quad \text{for} \ \ \pi\leq\theta<2\pi.\end{array}\right.$$
The spiking time of the neuron under this control is,
\begin{equation}
    T^M_{max}= \frac{2\pi+4\tan^{-1}\Big[z_dM/\sqrt{-z_d^2M^2+\omega^2}\Big]}{\sqrt{-z_d^2M^2+\omega^2}}, \label{eq:max_T_for_M<w/s}
\end{equation}
and the longest spiking time feasible with control $I^*$ is given by
\begin{equation}
    \label{eq:max_time_UBC}
    T^{I^*}_{max}=\int_0^{2\pi}\frac{1}{\sqrt{\omega^2+z_dM(z_dM-2\omega)\sin^2(\theta)}}.
\end{equation}
Then, by similar analysis for Case I, any spiking time $T\in[T^M_{min},T^{I^*}_{min})$ for a given $M<\w/z_d$ can be achieved with the minimum-power control $I^*_1$ as given in \eqref{eq:I1*}, any $T\in[T^{I^*}_{min},T^{I^*}_{max}]$ can be achieved with minimum power by $I^*$ in \eqref{eq:I*}, and moreover any $T\in(T^{I^*}_{max},T^M_{max}]$ can be obtained by switching between $I^*$ and $I^*_{max}$. The corresponding switching angles are $\theta_5=\sin^{-1}[2M\omega/(z_dM^2+z_d\omega \lambda_0)],\theta_6=\pi-\theta_5,\theta_7=\pi+\theta_5$ and $\theta_8=2\pi-\theta_5$, and the minimum-power optimal control for $T\in(T^{I^*}_{max},T^M_{max}]$ is characterized by
\begin{align*}
    \label{eq:I2*}
    I^*_2=\left\{\begin{array}{ll} I^* & \ 0\leq \theta < \theta_5 \\ -M & \ \theta_5\leq\theta\leq\theta_6 \\ I^* & \ \theta_6 < \theta < \theta_7 \\ M & \ \theta_7\leq\theta\leq\theta_8\\ I^* & \ \theta_8 < \theta \leq 2\pi.
    \end{array}\right.
\end{align*}
The $\lambda_0$ resulting in the optimal trajectory by $I_2^*$ can be calculated according to the given $T\in(T^{I^*}_{max},T^M_{max}]$ via the relation
\begin{equation*}
    T=\int_0^{\theta_5}{\frac{4}{\sqrt{\scriptstyle{\omega^2-\omega \lambda_0z_d^2\sin^2\theta}}}}d\theta+\int_{\theta_5}^{\frac{\pi}{2}}{\frac{4}{\omega-z_dM\sin\theta}}d\theta.
\end{equation*}

% =================== Figure 6 ====================
\begin{figure}[t]
\subfigure[]{\includegraphics[scale=0.105]{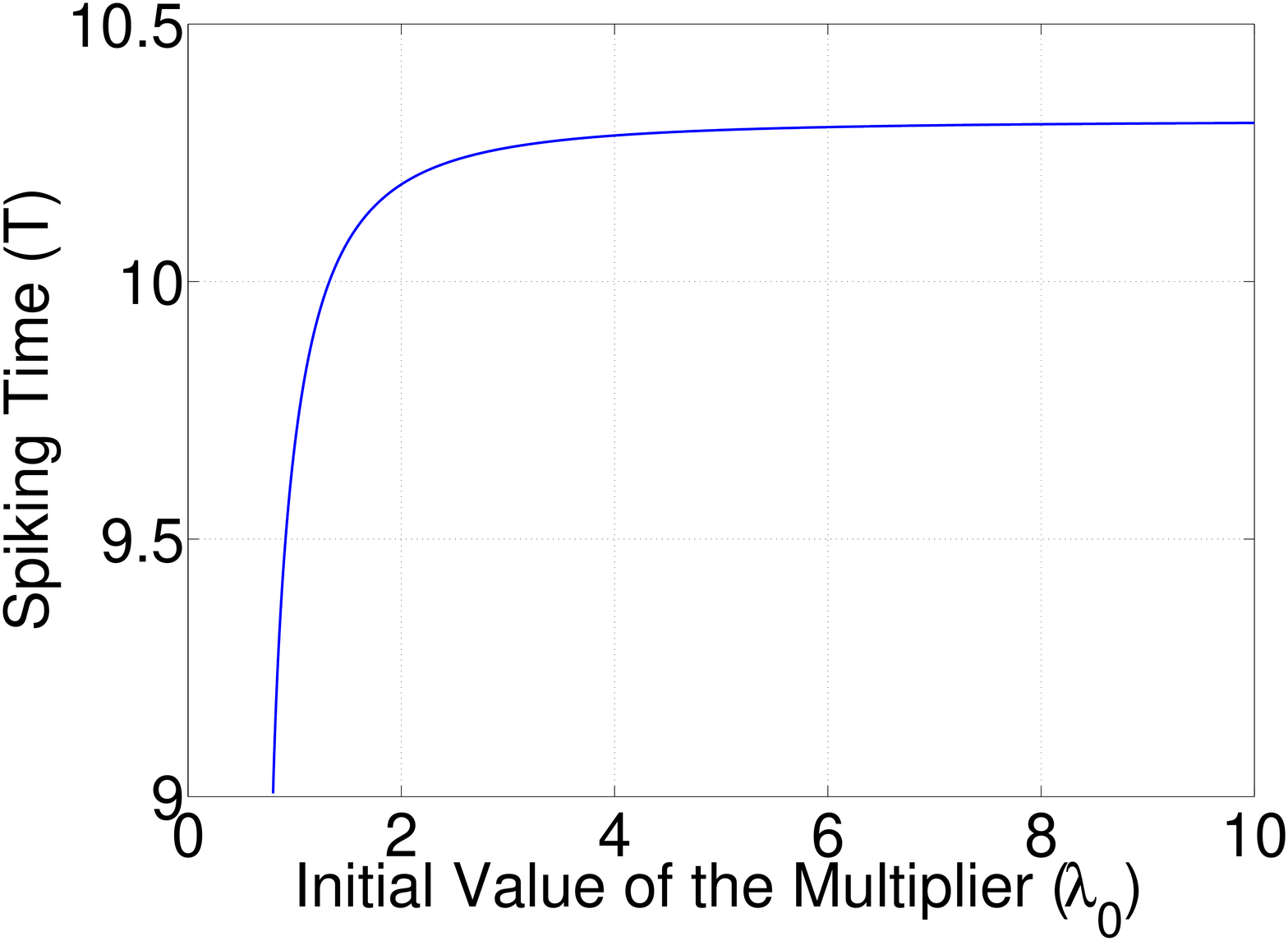}\label{fig:T_vs_lambda0_sine_for_I2}}
\subfigure[]{\includegraphics[scale=0.105]{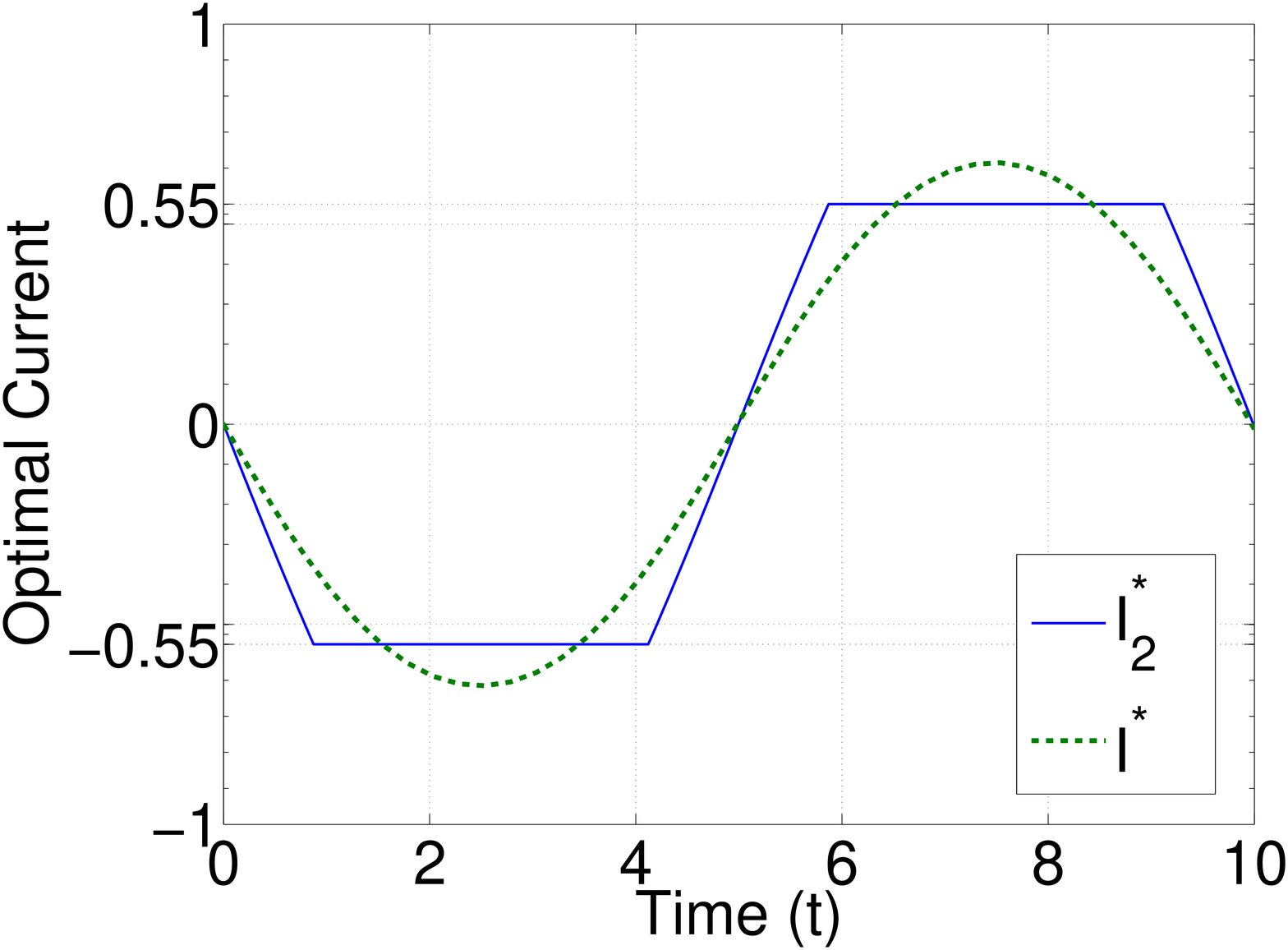}\label{fig:I2_and_I_sine}}
\caption{\subref{fig:T_vs_lambda0_sine_for_I2} Variation of the spiking time $T\in(T^{I^*}_{max},T^M_{max}]$ with the initial value of the multiplier, $\lambda_0$, for sinusoidal PRC model when $M=0.55$. \subref{fig:I2_and_I_sine} Minimum-power controls with ($M=0.55$) and without a constraint on the control amplitude for sinusoidal PRC model with $T=10$, $z_d=1$, and $\w=1$.}
	\label{fig:bounded_case2_sine}
\end{figure}

FIG. \ref{fig:T_vs_lambda0_sine_for_I2} shows the relation between $\lambda_0$ and $T$ by $I^*_2$ for $M=0.55$, $z_d=1$, and $\omega=1$. From \eqref{eq:max_T_for_M<w/s} the maximum possible spiking time with $M=0.55$ is $T^M_{max}=10.312$ and from \eqref{eq:max_time_UBC} the maximum spiking time feasible by $I^*$ is $I^{I^*}_{max}=9.006$. Therefore, in this example, any desired spiking time $T\in(9.006,10.312]$ can be obtained with minimum power by the use of $I^*_2$. FIG. \ref{fig:I2_and_I_sine} illustrates the bounded and unbounded optimal controls that spike the neuron at $T=10$, where $I^*$ is the minimum-power stimulus when the control amplitude is not limited and $I^*_2$ is the minimum-power stimulus when $M=0.55$. $I^*$ drives the neuron from $\t(0)=0$ to $\t(10)=2\pi$ with 2.193 units of power whereas $I^*_2$ requires 2.327 units.

\begin{figure}[ht]
\subfigure[]{\includegraphics[scale=0.42]{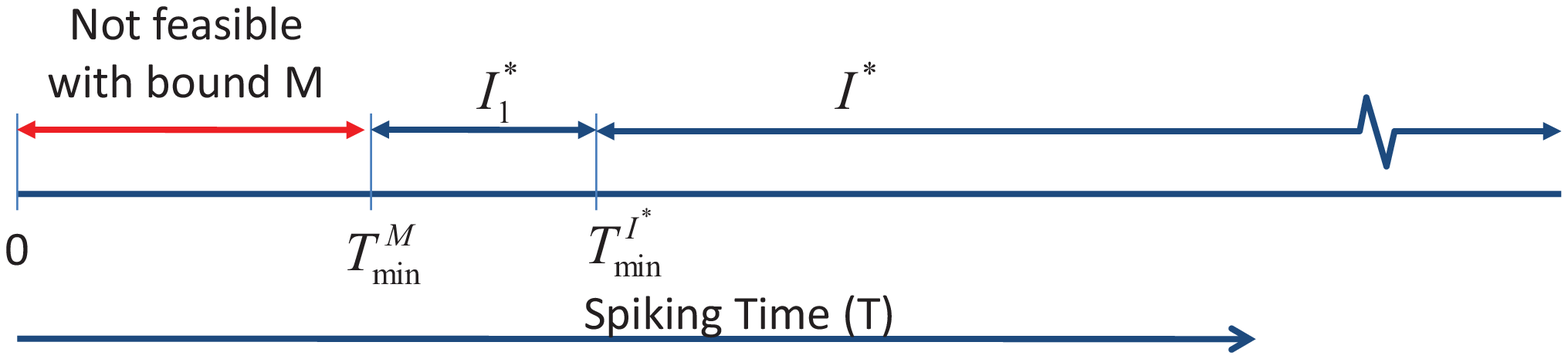}\label{fig:spiking_time_line_case1}}\\
\subfigure[]{\includegraphics[scale=0.42]{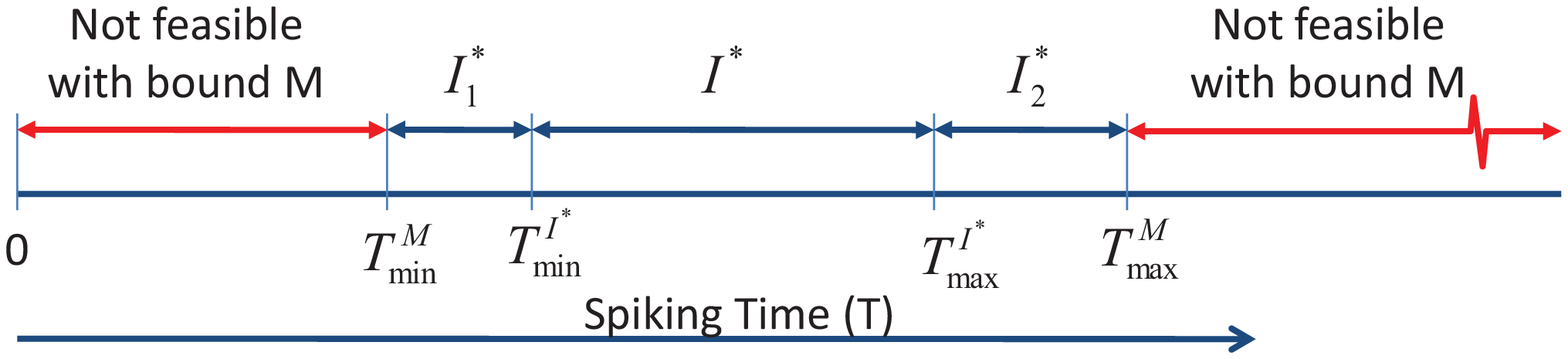}\label{fig:spiking_time_line_case2}}
\caption{A summery of the optimal control strategies for the sinusoidal PRC neuron for \subref{fig:spiking_time_line_case1} $M\geq\omega/z_d$, \subref{fig:spiking_time_line_case2} $M<\omega/z_d$. }
\label{fig:scenarios}
\end{figure}

A summary of the optimal (minimum-power) spiking scenarios for a prescribed spiking time of the neuron governed by the sinusoidal phase model \eqref{eq:sin_model} is illustrated in FIG. \ref{fig:scenarios}.

% ====================== Figure 7 ========================

%%%%%%%%%%%%%%%%%%%%%%%%%%%%%%%%%%%%%%%%%%%%%%%%%%%%%%%%%%%%%%%%%%%%%%%%%%%%%%%%%
\subsection{SNIPER PRC and Theta Neuron Phase Models}
\label{sec:sniper_prc}
We now consider the SNIPER PRC model in which $f(\theta)=\omega$ and $Z(\theta)=z_d(1-\cos\theta)$, where $z_d>0$ and $\omega>0$. That is,
\begin{equation}
     \label{eq:model_sniper}
     \dot{\theta}=\omega+z_d(1-\cos\theta)I(t).
\end{equation}
The minimum-power stimuli for spiking neurons modeled by this phase model can be easily derived with analogous analysis described previously in \ref{sec:unbounded_control_sine} and \ref{sec:bounded_control_sine} for the sinusoidal PRC phase model.

%%%%%%%%%%%%%%%%%%%%%%%%%%%%%%%%%%%%%%%%%%%%%%%%%%%%%%%%%%%%%%%%%%%%%%%%%%%%%%%%%
\subsubsection{Spiking Neurons with Unbounded Control}
\label{sec:unbounded_control_sniper}
Employing the maximum principle as in \ref{sec:unbounded_control_sine}, the minimum-power stimulus that spikes the SNIPER neuron at a desired time $T$ can be derived and given by
\begin{equation}
    \label{eq:I*sniper}
    I^*=\frac{-\omega+\sqrt{\omega^2-\omega \lambda_0 z_d^2(1-\cos\theta)^2}}{z_d(1-\cos\theta)},
\end{equation}
where $\lambda_0$ corresponding to the optimal trajectory is determined through the integral relation with $T$,
\begin{equation*}
    T=\int_0^{2\pi}{\frac{1}{\sqrt{\omega^2-\omega \lambda_0z_d^2(1-\cos \theta)^2}}}d\theta.
\end{equation*}

%============================ Figure8 =========================
\begin{figure}[t]
\subfigure[]{\includegraphics[scale=0.1]{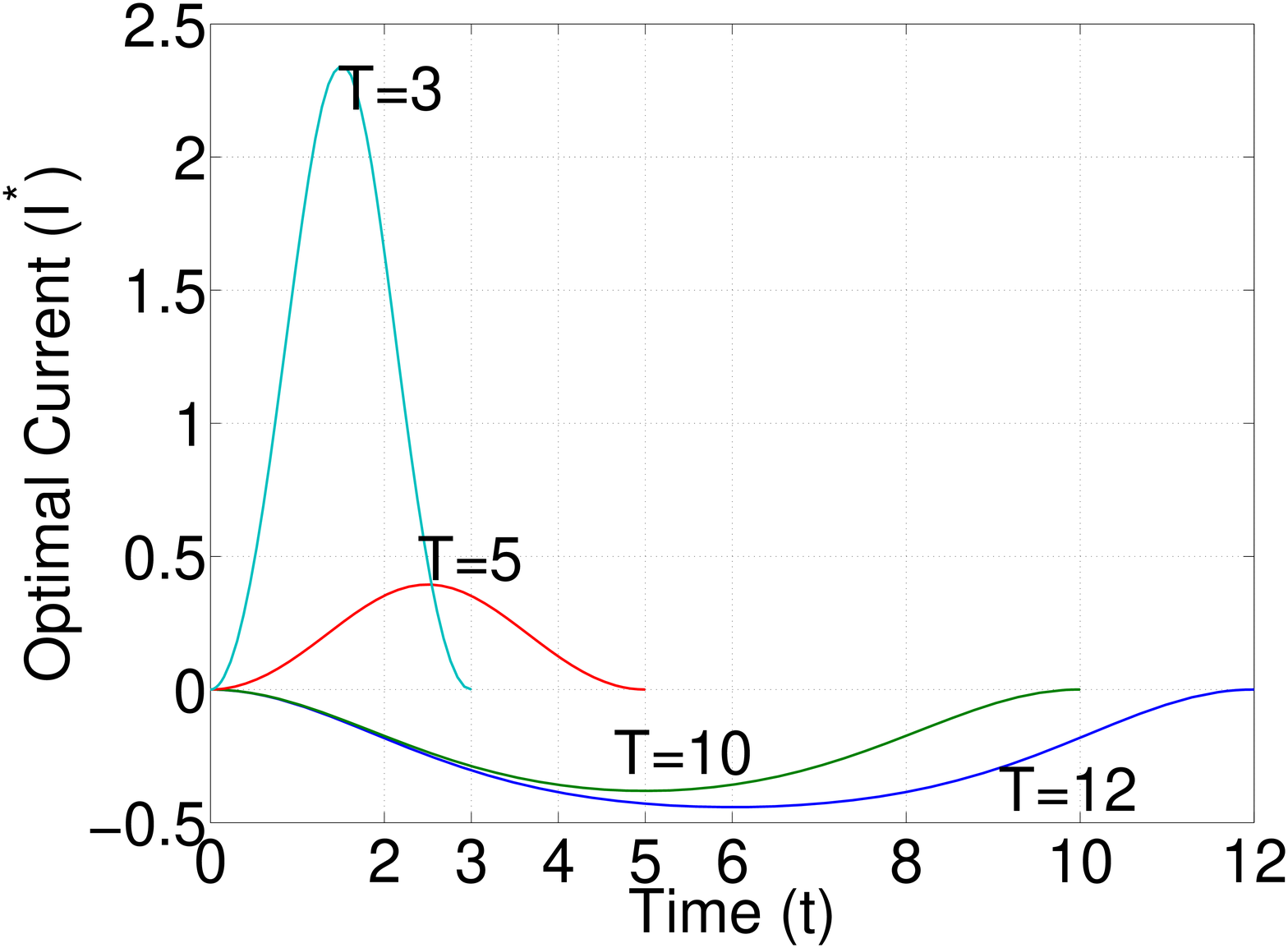}\label{fig:control_vs_time_sniper}}
\subfigure[]{\includegraphics[scale=0.1]{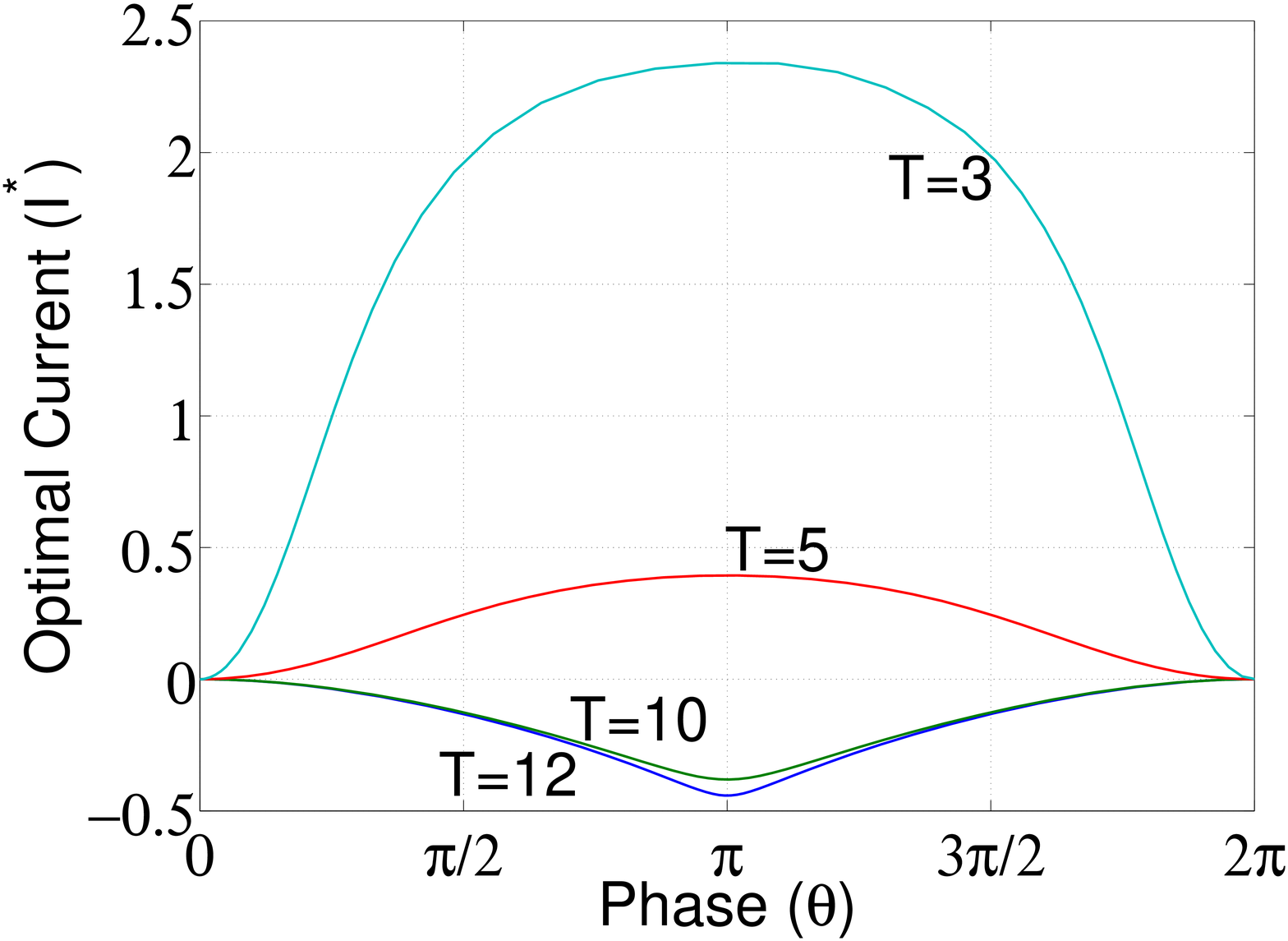}\label{fig:control_vs_theta_sniper}}\\ \subfigure[]{\includegraphics[scale=0.1]{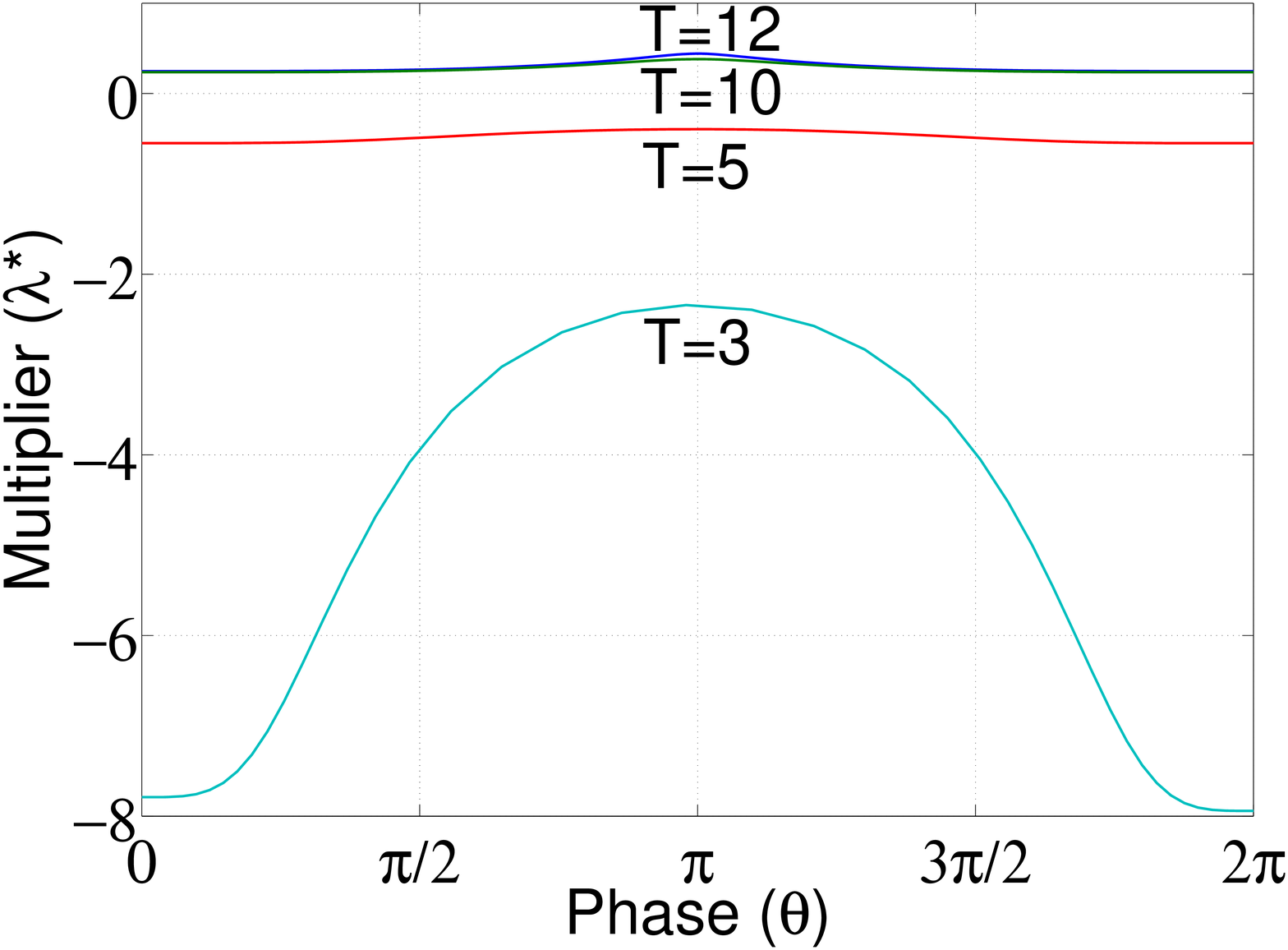}\label{fig:lambda_vs_theta_sniper}}
\subfigure[]{\includegraphics[scale=0.1]{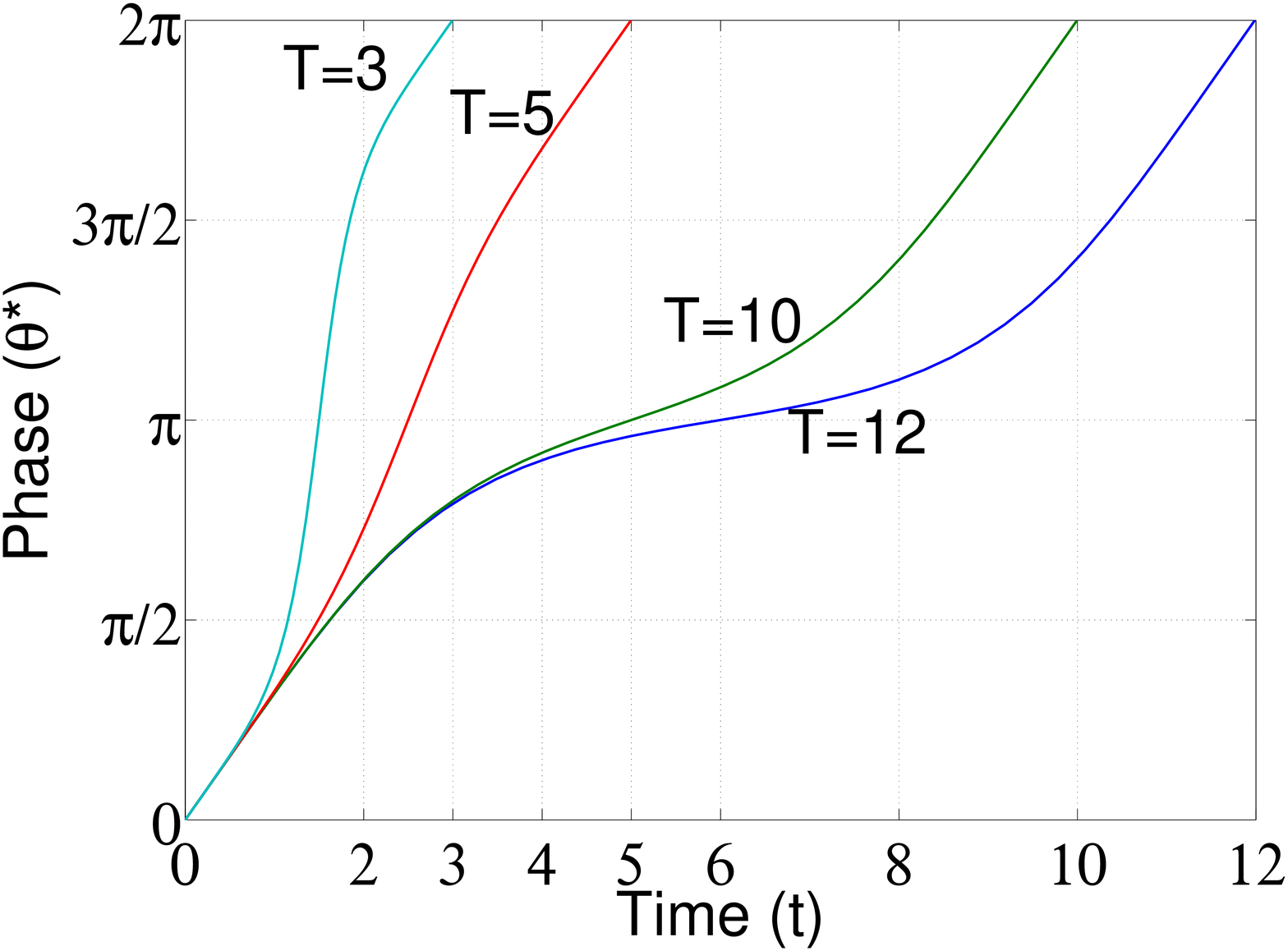}\label{fig:theta_vs_time_sniper}}
	  \caption{Optimal solutions for various spiking times $T=3,5,10,12$ for Sniper PRC model with $z_d=1$ and $\w=1$. \subref{fig:control_vs_time_sniper} The minimum-power control $I^*$. \subref{fig:control_vs_theta_sniper} Variation of $I^*$ with phase $\theta$. \subref{fig:lambda_vs_theta_sniper} Variation of the optimal multiplier, $\lambda^*$, with $\theta$. \subref{fig:theta_vs_time_sniper} Optimal phase trajectories following $I^*$.}
	\label{fig:trajectories_sniper}
\end{figure}

The minimum-power stimuli $I^*$ plotted with respect to time and phase for various spiking times $T=3,5,10,12$ with parameter values $z_d$=1 and $\omega=1$ are illustrated in FIG. \ref{fig:control_vs_time_sniper} and \ref{fig:control_vs_theta_sniper}, respectively. The corresponding optimal trajectories of $\lambda(\t)$ and $\t(t)$ for these spiking times are displayed in FIG. \ref{fig:lambda_vs_theta_sniper} and \ref{fig:theta_vs_time_sniper}.

%%%%%%%%%%%%%%%%%%%%%%%%%%%%%%%%%%%%%%%%%%%%%%%%%%%%%%%%%%
\subsubsection{Spiking Neurons with Bounded Control}
\label{sec:bounded_control_sniper}
When the amplitude of the available stimulus is limited, i.e., $|I(t)|\leq M$, the control that achieves the shortest spiking time for the SNIPER neuron modeled in \eqref{eq:model_sniper} is given by $I^*_{Tmin}=M>0$ for $0\leq\theta\leq 2\pi$, since $1-\cos\theta\geq 0$ for all $\theta\in [0,2\pi]$,. As a result, the shortest possible spiking time with this control is $T^M_{min}=2\pi/\sqrt{\omega^2+2z_d\omega M}$. Also, the shortest spiking time achieved by the control $I^*$ in \eqref{eq:I*sniper} given the bound $M$ is given by
\begin{equation}
    \label{eq:min_time_UBC_sniper}
    T^{I^*}_{min}=\int_0^{2\pi}\frac{1}{\sqrt{\omega^2+z_dM(z_dM+\omega)(1-\cos\theta)^2}}.
\end{equation}
Similar to the sinusoidal PRC case, the longest possible spiking time of the neuron varies with the control bound $M$. If $M\geq\omega/(2z_d)$, an arbitrarily large spiking time is achievable, however, if $M<\omega/(2z_d)$ there exists a maximum spiking time.

\begin{figure}[t]
\subfigure[]{\includegraphics[scale=0.1]{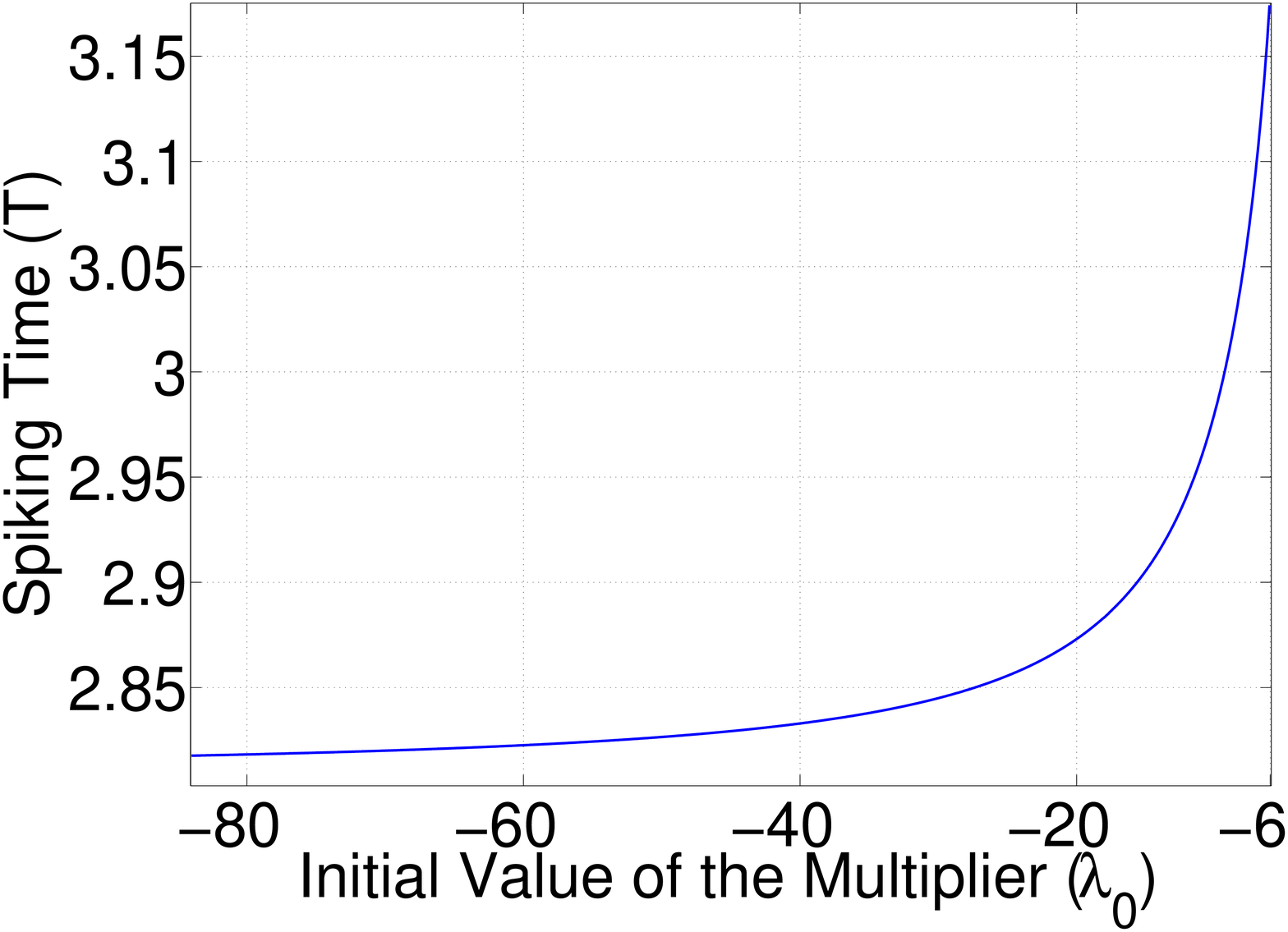}\label{fig:T_vs_lambda0_sniper_for_I1}}
\subfigure[]{\includegraphics[scale=0.1]{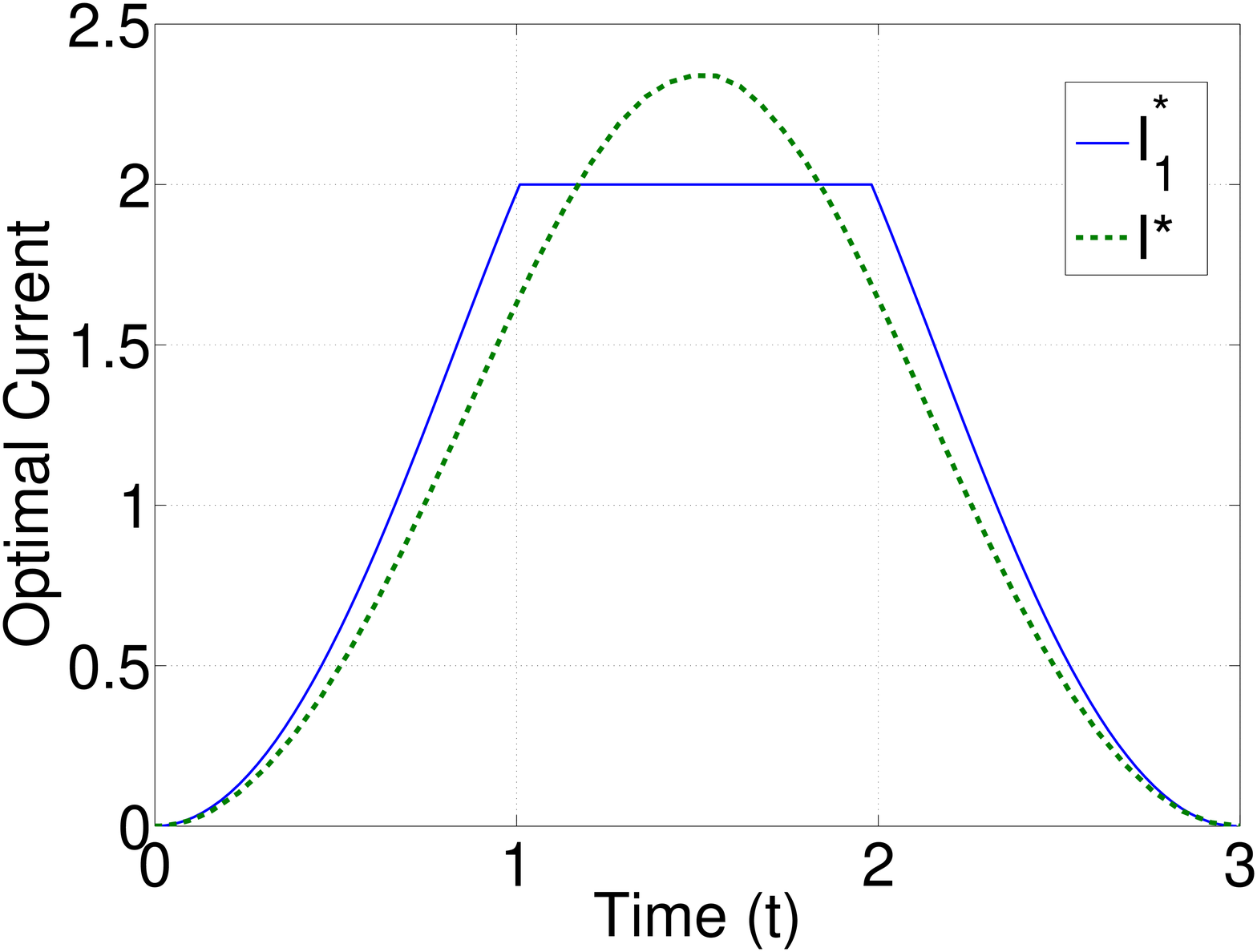}\label{fig:I1_and_I_sniper}}
\caption{\subref{fig:T_vs_lambda0_sniper_for_I1} Variation of the spiking time $T\in[T_{min},T'_{min})$ with the initial multiplier value, $\lambda_0$, for SNIPER PRC model when $M=2$. \subref{fig:I1_and_I_sniper} Minimum-power controls with ($M=2$) and without a constraint on the control amplitude for SNIPER PRC model with $T=3$, $z_d=1$, and $\omega=1$.}
	\label{fig:bounded_sniper_case1}
\end{figure}

% ----------------------------------------
\emph{Case I: $M\geq\omega/(2z_d)$}. Any spiking time $T\in[T^{I^*}_{min},\infty)$ is possible with control $I^*$ but a shorter spiking time $T\in[T^M_{min},T^{I^*}_{min})$ requires switching between $I^*$ and $I^*_{Tmin}$, which is characterized by two switchings,
\begin{align}
    \label{eq:I1*sniper}
    I^*_1=\left\{\begin{array}{ll} I^*,&\quad 0\leq \theta < \theta_1\\
    M,&\quad \theta_1\leq\theta\leq 2\pi-\theta_1 \\
    I^*,& \quad 2\pi-\theta_1 < \theta \leq 2\pi
\end{array}\right.
\end{align}
where $\theta_1=\cos^{-1}\left[1+2\omega M/(z_dM^2+z_d\omega \lambda_0)\right]$. The initial value $\lambda_0$ which results in the optimal trajectory is given by,
\begin{equation*}
T=\int_0^{\theta_1}{\frac{2}{\sqrt{\scriptstyle{\omega^2-\omega \lambda_0z_d^2(1-\cos\theta)^2}}}}d\theta+\int_{\theta_1}^{\pi}{\frac{2}{\scriptstyle{\omega+z_dM(1-\cos\theta)}}}d\theta.
\end{equation*}
FIG. \ref{fig:T_vs_lambda0_sniper_for_I1} illustrates the relation between $\lambda_0$ and $T\in[T^M_{min},T^{I^*}_{min})$ by $I_1^*$ for $M=2$, $z_d=1$, and $\omega=1$. In this case, the shortest feasible spiking time is $T^{M}_{min}=2.09$ and the shortest with the control $I^*$ is $T^{I^*}_{min}=3.18$. Any spiking time in the interval $(2.09,3.18]$ is achievable by $I^*_1$ in \eqref{eq:I1*sniper} with minimum-power. FIG. \ref{fig:I1_and_I_sniper} illustrates the unbounded and bounded, with $M=2$, optimal stimuli that fire the neuron at $T=3$ with minimum-power.

%========================= Figure9 ===========================

\begin{figure}[t]
\subfigure[]{\includegraphics[scale=0.1]{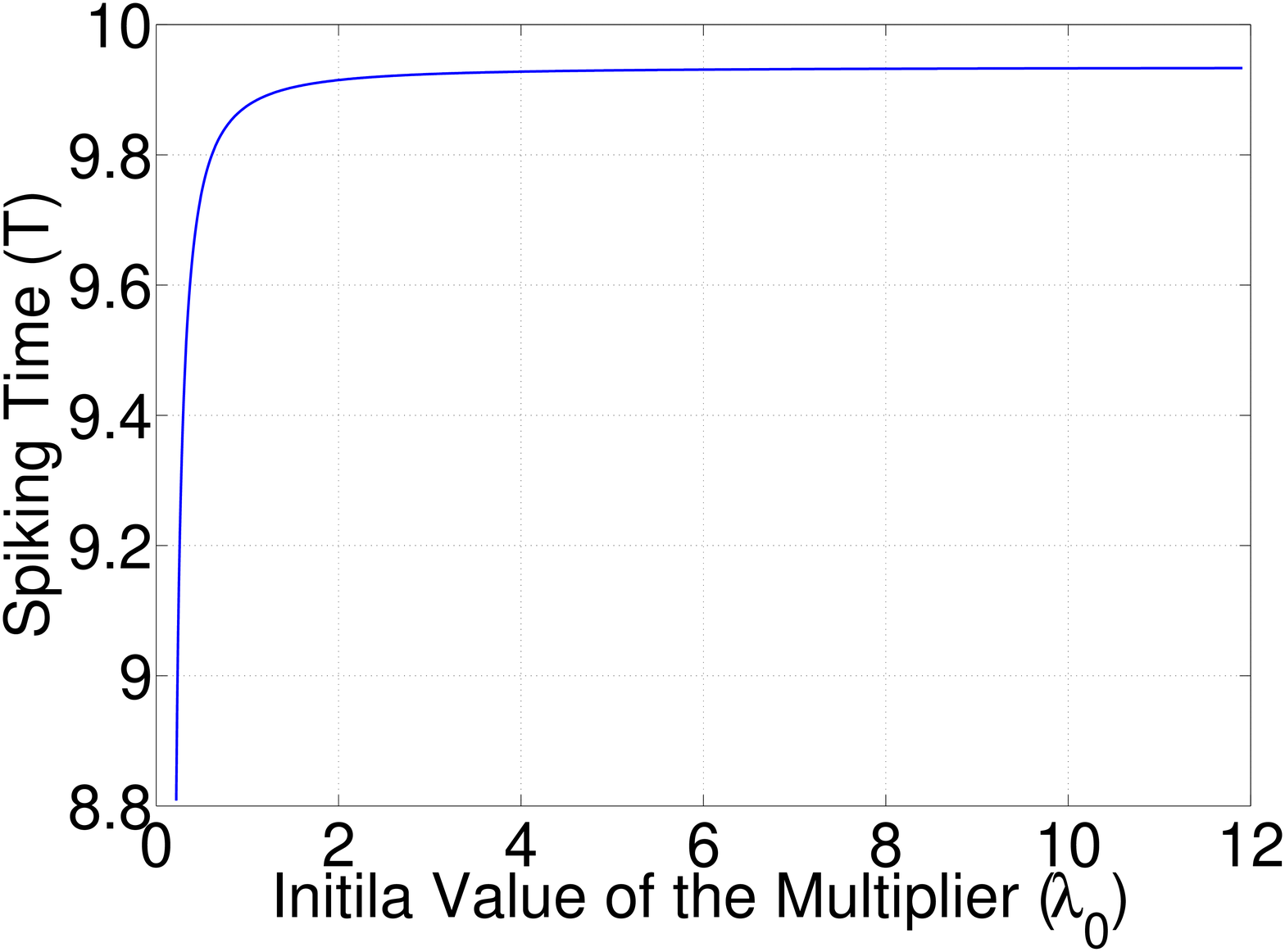}\label{fig:T_vs_lambda0_sniper_for_I2}}
\subfigure[]{\includegraphics[scale=0.1]{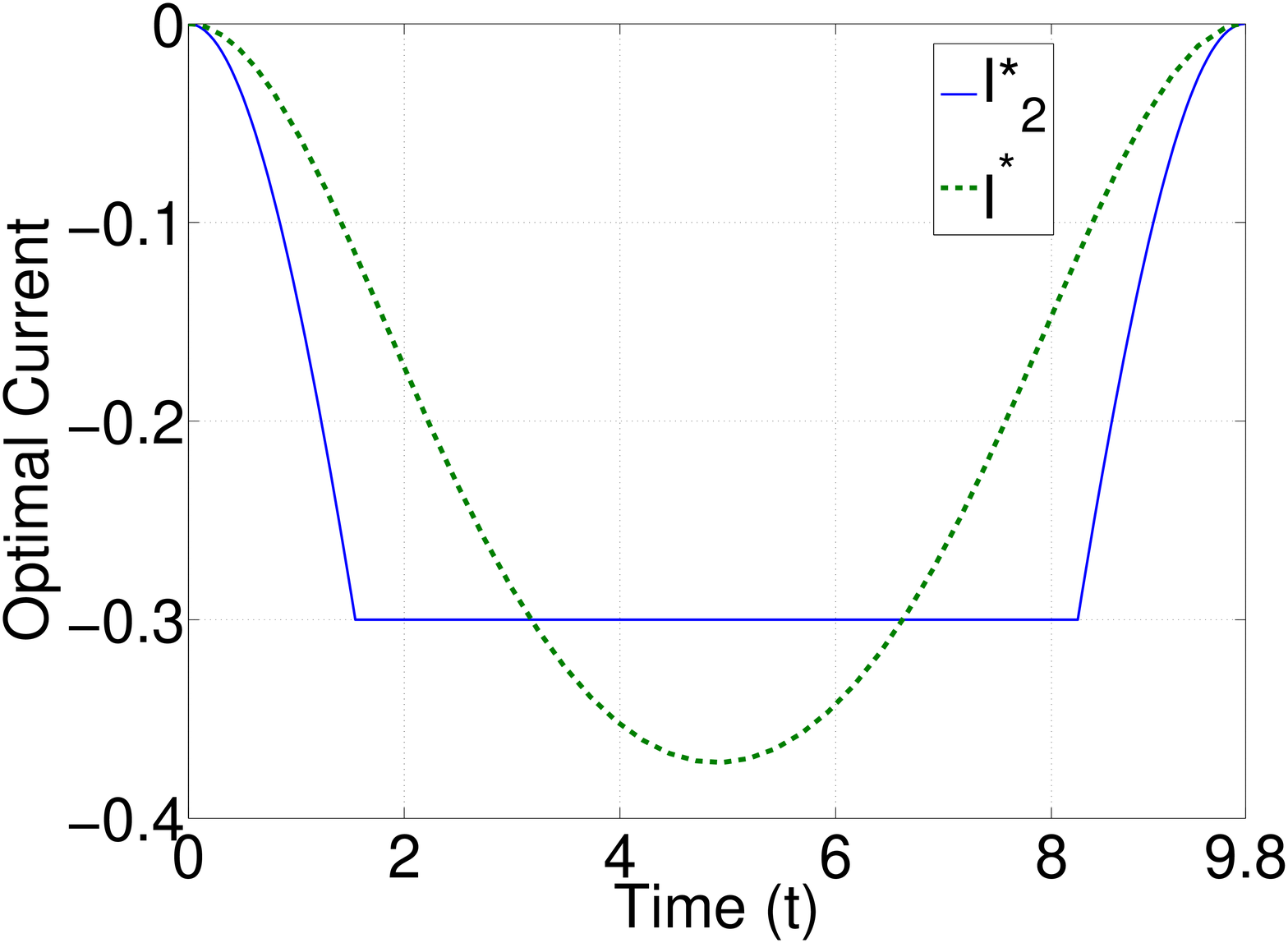}\label{fig:I2_and_I_sniper}}
 	  \caption{\subref{fig:T_vs_lambda0_sniper_for_I2} Variation of the spiking time $T\in(T'_{max},T_{max}]$ with the initial multiplier value, $\lambda_0$, for SNIPER PRC model when $M=0.3$. \subref{fig:I2_and_I_sniper} Minimum-power controls with ($M=0.3$) and with out a constraint on the control amplitude for SNIPER PRC model with $T=9.8$, $z_d=1$, and $\omega=1$.}
	\label{fig:bounded_sniper_case2}
\end{figure}

% -------------------------------------------
\emph{Case II: $M<\omega/(2z_d)$}.
In this case there exists a longest possible spiking time which is achieved by $I_{max}^*=-M$ for all $\theta\in [0,2\pi]$. The longest spiking time feasible with the control $I^*$ as in \eqref{eq:I*sniper} is given by
\begin{equation*}
    T^{I^*}_{max}=\int_0^{2\pi}\frac{1}{\sqrt{\omega^2+z_dM(z_dM-2\omega)(1-\cos\theta)^2}}.
\end{equation*}
Therefore, any spiking time $T\in[T^M_{min},T^{I^*}_{min})$ for a given $M<\w/z_d$ can be achieved with the minimum-power control $I^*_1$ as given in \eqref{eq:I1*sniper}, any $T\in[T^{I^*}_{min},T^{I^*}_{max}]$ can be achieved with minimum power by $I^*$ in \eqref{eq:I*sniper}, and moreover any $T\in(T^{I^*}_{max},T^M_{max}]$ can be obtained by switching between $I^*$ and $I^*_{max}$, that is,
\begin{align*}
    I^*_2=\left\{\begin{array}{ll} I^*,& \quad 0\leq \theta < \theta_2  \\
    -M,& \quad \theta_2\leq\theta\leq 2\pi-\theta_2 \\
    I^*,&\quad 2\pi-\theta_2 < \theta < 2\pi
\end{array}\right.
\end{align*}
where $\theta_2=\cos^{-1}\left[1-2\omega M/(z_dM^2+z_d\omega \lambda_0)\right]$. The $\lambda_0$ associated with the optimal trajectory is determined via the relation with the desired spiking time $T$,
\begin{equation*}
    T=\int_0^{\theta_1}{\frac{2}{\sqrt{\scriptstyle {\omega^2-\omega \lambda_0z_d^2(1-\cos \theta)^2}}}}d\theta+\int_{\theta_1}^{\pi}{\frac{2}{\scriptstyle {\omega-z_dM(1-\cos\theta)}}}d\theta.
\end{equation*}
FIG. \ref{fig:T_vs_lambda0_sniper_for_I2} illustrates the relation between $\lambda_0$ and $T\in(T^{I^*}_{max},T^M_{max}]$ by $I_2^*$ for $M=0.3$, $z_d=1$, and $\omega=1$. In this case, the longest feasible spiking time is $T^M_{max}=9.935$ and the longest with the control $I^*$ is $T^{I^*}_{max}=8.596$. The unbounded and bounded, with $M=0.3$, optimal stimuli that fire the neuron at $T=9.8$ with minimum-power are illustrated in FIG. \ref{fig:I2_and_I_sniper}.

A summary of the optimal (minimum-power) spiking scenarios for a prescribed spiking time of the neuron governed by the SNIPER PRC model in \eqref{eq:model_sniper} can be illustrated analogously to FIG. \ref{fig:spiking_time_line_case1} and \ref{fig:spiking_time_line_case2} for $M\geq \omega/(2z_d)$ and $M<\omega/(2z_d)$, respectively.

%============================ Figure10 ===============================

% ============================== Theta Neuron ======================================
The theta neuron phase model is described by the dynamical equation
\begin{equation}
    \label{eq:model_theta}
    \dot{\theta}=1+\cos\theta+z_d(1-\cos\theta)(I(t)+I_b),
\end{equation}
where $I_b>0$ is the baseline current and the natural frequency $\omega$ of the theta neuron is given by $2\sqrt{I_b}$. In fact, by a coordinate transformation, $\theta(\phi)=2\tan^{-1}\left[\sqrt{I_b}\tan\left((\phi-\pi)/2\right)\right]+\pi$, the theta neuron model can be transformed to a SNIPER PRC with $z_d=\w/2$ in \eqref{eq:model_sniper}. Therefore, all of the results developed for the SNIPER PRC can be directly applied to the theta neuron phase model.

%%%%%%%%%%%%%%%%%%%%%%%%%%%%%%%%%%%%%%%%%%%%%%%%%%%%%%%%%%%%%%%%%%%%%%%%%%%%%%%%%%%%%%%%%%%%%%%%%%%%
\section{Conclusion and Future Work\label{sec:conclusion}}
In this paper, we studied various phase-reduced models that describe the dynamics of a neuron system. We considered the design of minimum-power stimuli for spiking a neuron at a specified time instant and formulated this as an optimal control problem. We investigated both cases when the control amplitude is unbounded and bounded, for which we found analytic expressions of optimal feedback control laws. In particular for the bounded control case, we characterized the ranges of possible spiking periods in terms of the control bound. Moreover, minimum-power stimuli for steering any nonlinear oscillators of the form as in \eqref{eq:phasemodel} between desired states can be derived following the steps presented in this article. In addition, the charge-balanced constraint \cite{Nabi09} can be readily incorporated into this framework as well.

The optimal control of a single neuron system investigated in this work illustrates the fundamental limit of spiking a neuron with external stimuli and provides a benchmark structure that enables us to study optimal control of a spiking neural network with different individual oscillation frequencies. Our recent work \cite{Li_NOLCOS10} proved that simultaneous spiking of a network of neurons is possible; however, optimal control of such a spiking  neural network has not been studied. We finally note that although one-dimensional phase models are reasonably accurate to describe the dynamics of neurons, studying higher dimensional models such as that of Hodgkin-Huxley is essential for more accurate computation of optimal neural inputs.

\nocite{*}
\bibliography{SingleNeuron_PRE}% Produces the bibliography via BibTeX.
\end{document}